\documentclass[a4paper,9pt,parskip=half,english]{article}
\usepackage[top=60pt,bottom=60pt,left=78pt,right=76pt]{geometry}
\usepackage{babel}
\usepackage[T1]{fontenc}
\usepackage{lmodern,textcomp,latexsym}
\usepackage{amsthm,amssymb,amsmath,amsfonts,dsfont,mathrsfs,bbm,bm,mathtools}
\usepackage{graphicx,psfrag,color,rotate}
\usepackage{paralist,threeparttable}
\usepackage[footnotesize]{caption}
\usepackage{dsfont,booktabs,threeparttable}
\usepackage{longtable,booktabs}
\usepackage{cite}
\usepackage{tikz}
\usepackage{pgfplotstable}
\usepackage{pgfplots}
\usepackage{nicematrix}
\usepackage[title]{appendix}
\usetikzlibrary{calc,spy,external}
\usepackage{caption}
\usepackage{subcaption}
\usepackage{float}

\usepackage{titlesec}

\usepackage{booktabs}

\usepackage{gettitlestring}
\usepackage{hyperref}

\titleformat*{\section}{\large\bfseries}
\titleformat*{\subsection}{\large\bfseries}
\titleformat*{\subsubsection}{\large\bfseries}
\titleformat*{\paragraph}{\large\bfseries}
\titleformat*{\subparagraph}{\large\bfseries}

\setdefaultenum{1)}{\theenumi.1)}{}{}

\theoremstyle{plain}
\newtheorem{theorem}{Theorem}
\newtheorem{proposition}{Proposition}
\newtheorem{remark}[theorem]{Remark}
\theoremstyle{definition}
\newtheorem{definition}{Definition}

\newcommand{\T}{^{\mathop{\mathrm{T}}}}
\newcommand{\negTranspose}{^{\mathop{\mathrm{-T}}}}

\newcommand{\mR}{\mathbb{R}}
\newcommand{\mRnn}{\mathbb{R}^{n \times n}}
\newcommand{\mRNN}{\mathbb{R}^{N \times N}}

\newcommand{\vxi}{\bm{\xi}}

\newcommand{\vph}{\bm{\varphi}}

\newcommand{\vPh}{\mathbf \Phi}

\newcommand{\vb}{\mathbf b}
\newcommand{\vc}{\mathbf c}

\newcommand{\ve}{\mathbf e}
\newcommand{\vf}{\mathbf f}
\newcommand{\vg}{\mathbf g}

\newcommand{\vp}{\mathbf p}
\newcommand{\vq}{\mathbf q}

\newcommand{\vv}{\mathbf v}
\newcommand{\vw}{\mathbf w}
\newcommand{\vx}{\mathbf x}
\newcommand{\vy}{\mathbf y}
\newcommand{\vz}{\mathbf z}

\newcommand{\vA}{\mathbf A}
\newcommand{\vB}{\mathbf B}
\newcommand{\vC}{\mathbf C}

\newcommand{\vF}{\mathbf F}

\newcommand{\vH}{\mathbf H}
\newcommand{\vI}{\mathbf I}

\newcommand{\vK}{\mathbf K}

\newcommand{\vM}{\mathbf M}

\newcommand{\vV}{\mathbf V}

\newcommand{\vX}{\mathbf X}

\def\addlegendimage{\csname pgfplots@addlegendimage\endcsname}

\pgfplotsset{
	cycle list={%
		{draw=black,mark=star,solid},
		{draw=black, mark=square,solid},
		{draw=black,mark=+,solid},
		{black,mark=o},}}

\begin{document}
	\begin{center}
		\textbf{\large Extended invariant cones as Nonlinear Normal Modes of inhomogeneous piecewise linear systems }
		
		\qquad

		\renewcommand{\thefootnote}{\fnsymbol{footnote}}
		A.Yassine Karoui$^{1}$\footnote{karoui@inm.uni-stuttgart.de},
		and
		Remco I. Leine$^{1}$\footnote{leine@inm.uni-stuttgart.de}
		
		\qquad
		
		$^{1}$Institute for Nonlinear Mechanics
		
		University of Stuttgart
		
		Pfaffenwaldring 9, 70569 Stuttgart, Germany

	\end{center}
	\vspace{5pt}
	\hrule
	\vspace{8pt}
	\textbf{Abstract}: The aim of this paper is to explore the relationship between invariant cones and nonlinear normal modes in piecewise linear mechanical systems. As a key result, we extend the invariant cone concept, originally established for homogeneous piecewise linear systems, to a class of inhomogeneous continuous piecewise linear systems.  The inhomogeneous terms can be constant and/or time-dependent, modeling nonsmooth mechanical systems with a clearance gap and external harmonic forcing, respectively. Using an augmented state vector, a modified invariant cone problem is formulated and solved to compute the nonlinear normal modes, understood as periodic solutions of the underlying conservative dynamics. An important contribution is that invariant cones of the underlying homogeneous system can be regarded as a singularity in the theory of nonlinear normal modes of continuous piecewise linear systems. In addition, we use a similar methodology to take external harmonic forcing into account. We illustrate our approach using numerical examples of mechanical oscillators with a unilateral elastic contact. The resulting backbone curves and frequency response diagrams are compared to the results obtained using the shooting method and brute force time integration.  \\[8pt]
	\textit{Keywords}: inhomogeneous piecewise linear systems, invariant cones, nonlinear normal modes, backbone curves, frequency response curves, shooting method.  
	\vspace{8pt}
	\hrule
\setcounter{footnote}{0}
\renewcommand{\thefootnote}{\roman{footnote}}
\section{Introduction}\label{section:intro}
The objective of this paper is to relate the invariant cones of a class of continuous piecewise linear systems to the framework of nonlinear normal modes. We propose an extended formulation of the invariant cone problem taking inhomogeneous terms into account, with which nonlinear normal modes of piecewise linear systems with and without clearance gaps, damping and external harmonic forcing can be computed. Thereby, the connection between the invariant cone problem and the concept of nonlinear normal modes is unveiled. 

The concept of nonlinear normal modes (NNMs) is fundamental to understanding the behavior of nonlinear dynamical systems and has been extensively explored for a wide range of system classes. Originally, the term "nonlinear mode" dates back to the seminal work of Rosenberg in the sixties \cite{Rosenberg1960, Rosenberg1961, Rosenberg1962, Rosenberg1964, Rosenberg1964a, Rosenberg1966}, who defined NNMs as synchronous periodic oscillations of conservative nonlinear systems, seen as a natural extension of linear modes (LNMs). This definition requires that all material points of the system reach their extreme values and pass through zero simultaneously, allowing all displacements to be expressed in terms of a single reference displacement with no velocity dependence. In this context, the modal lines in the configuration space are generally curved, their shape depending on the total energy in the system due to the nonlinear relationship between coordinates during periodic motion. However, this definition does not extend to nonconservative systems. NNMs were further studied in the seventies by Rand \cite{Rand1971,Rand1971a,Rand1974} and Manevitch and Mikhlin \cite{Manevich1972}.
In the nineties, the NNM research gained momentum through the efforts of Vakakis et al. \cite{Caughey1990,Vakakis1992,Vakakis1996,Vakakis1997} and Shaw and Pierre \cite{Shaw1991,Shaw1993,Shaw1994}. In \cite{Shaw1993}, Shaw and Pierre introduced a new definition of NNMs by relaxing the synchronicity condition and applying the center manifold approach to vibrating systems. They defined NNMs as two-dimensional manifolds tangent to the linear subspaces in phase space. These manifolds are invariant under the flow and are parameterized using a pair of state variables, consisting of both displacement and velocity, while assuming a functional relationship of the remaining variables to the chosen pair. This definition allowed to broaden the spectrum of nonlinear systems that can be studied using NNMs by including effects such as damping \cite{Shaw1993} and later external forcing \cite{Jiang2005}.

Even though smooth nonlinear systems were initially at the center of focus, many studies have been conducted on nonsmooth mechanical systems of piecewise nature with the purpose of understanding and computing their NNMs. One of the first and most important efforts in this direction was the work by Zuo and Curnier \cite{Zuo1994}, who extended Rosenberg´s definition to conservative unforced mechanical systems with piecewise linear (hereafter called PWL) restoring forces. They studied a class of PWL systems with the nonlinearity located at the origin, being the equilibrium of the system, i.e. mechanical oscillators with a one-sided spring with zero clearance. For such systems, they found that the NNM frequencies are energy-independent, as is the case in linear systems. Later on, Chen and Shaw \cite{Chen1996} proposed a method to construct NNMs of unforced conservative PWL systems, by applying the invariant manifold approach from \cite{Shaw1993} to a nonsmooth system with a nonvanishing clearance gap. This approach was then extended in \cite{Jiang2004} to account for large amplitude vibrations. Therein, the authors propose a numerical technique to compute NNMs for piecewise linear autonomous systems based on a Galerkin method. Other investigations, such as \cite{Chati1997, Butcher1999} and \cite{Uspensky2014}, explored similar classes of unforced and undamped nonsmooth systems, whereas external forcing was considered in \cite{Butcher2007,Saito2011,Uspensky2013} and \cite{Uspensky2019}. Regarding the numerical computation of NNMs of PWL systems, several state-of-the-art methods have been extended from the smooth framework to take nonsmooth nonlinearities into account. For instance, the harmonic balance method combined with the asymptotic numerical method can be used to compute branches of NNMs for PWL systems, provided that the nonsmoothness is regularized, as shown in \cite{Cochelin2014}, \cite{Karkar2013} and \cite{Moussi2015}. A simpler and more popular approach is the shooting method coupled with an arclength-continuation algorithm to trace NNMs branches. Originally introduced for smooth systems, the shooting method is a periodic-solution solver in time domain based on the Newton method. The technique performs an iterative improvement of initial estimates for the nonlinear mode shapes and frequencies by enforcing the free response of the system to satisfy a periodicity condition complemented by arbitrary constraints. Using a slightly modified procedure enables the shooting method to handle nonsmooth nonlinearities without requiring a regularization. An example is provided in the study by Attar et al. \cite{Attar2017}, where thin-walled structures with unilateral contact forces were considered. The computation of NNMs branches using the shooting method combined with a pseudo-arclength continuation technique is described in detail for a general class of nonlinear systems in \cite{Kerschen2014}.

The common aim of all the mentioned investigations consisted in finding invariant sets of continuous nonsmooth systems, a task which proved more challenging than in the smooth framework due to the lack of differentiability. In piecewise linear systems with a non-vanishing gap, which we refer to as inhomogeneous, the equilibrium lies in a region of the state-space characterized by locally linear behavior and the nonlinearity is triggered at a specific amplitude. In contrast, in the homogeneous case, the equilibrium is located on the switching manifold which makes a local linearization of the system in the form $\dot{x} = A x$ impossible, thereby excluding the possibility to characterize the local behavior using eigenvectors and eigenvalues. For homogeneous systems, invariant cones were identified as a powerful tool to investigate bifurcations and stability. Originally, the concept of an "invariant cone" was introduced by Carmona et al. \cite{carmona2005bifurcation}, who studied a three-dimensional continuous piecewise linear system with two subregions, hereafter called 2CPL$_3$ systems. In \cite{carmona2005bifurcation}, the authors identified invariant cones as sets in the state-space characterized by invariant half-lines on a Poincaré section. The existence of an unstable invariant cone explains the perhaps counter-intuitive behavior of a piecewise linear system composed of two stable subsystems that combine to form an unstable system. This behavior is not exclusive for low dimensional academic example systems and was observed in the context of nonsmooth mechanics by Rossi et al. \cite{Rossi2023}, who recently observed a flutter-like instability in a mechanical system with follower forces and explained its unstable behavior through the existence of an unstable invariant cone. This important investigation was preceded by the work of Küpper and coworkers \cite{Kuepper2008, Kuepper2011, Hosham2016}, who extended Carmona´s idea to higher dimensional piecewise linear systems, both continuous and discontinuous. 

To the best of our knowledge, no connection between the concepts of invariant cones and nonlinear normal modes has been established, despite the apparent similarity between the two. This paper aims to bridge this gap by extending the invariant cone formulation to a broader class of inhomogeneous PWL systems and using this modified formulation to compute NNMs. For this, NNMs are regarded as periodic solutions of the underlying conservative dynamics, represented as frequency-energy plots often called "backbone curves" or of the damped and forced dynamics in the sense of \cite{Jiang2005}. An augmented state vector is introduced to transform the inhomogeneous PWL system including constant and time-dependent forcing terms into a homogeneous form. This augmented dynamics is used to formulate a modified version of the invariant cone problem following its definition in \cite{Kuepper2008}, whose solution tracks the nonlinear regions of the NNM branches. Special attention is given to the inhomogeneity of the systems under investigation. In fact, in its original form, the concept of invariant cones is based on the positive homogeneity property, which leads to constant time intervals in which orbits oscillate in one of the subregions. This property is destroyed by the presence of an inhomogeneous term in the dynamics, leading to amplitude-dependent return times. This, in turn, is expected from a mechanical point of view, since nonlinear normal modes of inhomogeneous PWL systems are periodic orbits with amplitude-dependent frequencies. The interplay between the inhomogeneity in our modified invariant cone approach and the characteristics of the NNMs of inhomogeneous PWL systems is addressed in detail in this work. 

The paper is organized as follows. We introduce the class of continuous PWL systems of interest in Section 2. The theoretical background of invariant cones is summarized in Section 3. Section 4 extends the invariant cone concept to inhomogeneous autonomous PWL systems, where the inhomogeneous term is constant and is due to a non-vanishing clearance gap. Therein, a modified version of an invariant cone problem is solved to trace the nonlinear regions of the backbone curves and comparisons with the autonomous shooting method are drawn. Finally, in Section 5, damping and external harmonic forcing are considered in systems with or without a clearance gap. Here, another modified invariant cone problem is proposed, whose numerical solution coupled with a pseudo-arclength continuation technique allows to compute nonlinear frequency-response curves. We conclude by summarizing
our contributions and discussing future work.

\section{Continuous PWL systems}\label{section:background}
\begin{figure}
	\centering	
	\includegraphics[scale=1]{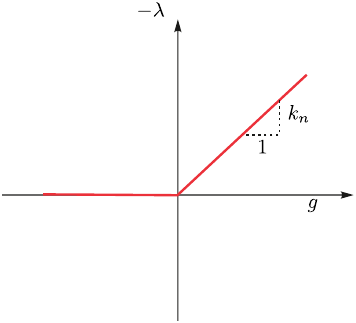}	
	\caption{Characteristic of the unilateral elastic contact law.}
	\label{fig_unilateral_elastic_force_law}
\end{figure}

To begin with, we formally introduce the family of continuous piecewise linear systems to be studied here, based on the work by Carmona \cite{Carmona2002}.
\begin{definition}
	A differential equation \begin{equation}
		\dot{\vx} = \vF(\vx), \nonumber
	\end{equation}
	with $\vx\in \mR^n$ is said to be a 2CPL$_n$ system if the right-hand side $\vF: \mR^n \mapsto \mR^n$ is continuous and of the form
	\begin{equation}
		\dot{\vx} = \vF(\vx) = \begin{cases}
			\vA^- \vx + \vb^-, \quad \mbox{if } \vv\T \vx < \delta \\ 
			\vA^+ \vx + \vb^+, \quad \mbox{if } \vv\T \vx \geq \delta 
		\end{cases},
	\end{equation}
	where $\vA^-, \vA^+ \in \mRnn$, $\vb^\pm \in \mR^n$, $\delta \in \mR$ and $\vv \in \mR^n$ is a non-vanishing vector. The continuity of $\vF$ requires that for all $\vx^*$ such that $\vv\T \vx^* = \delta $ it holds that $\vA^-\vx^* + \vb^- = \vA^+ \vx^* + \vb^+ $, therefore putting conditions on the matrices $\vA^\pm$ and the vectors $\vb^\pm$.  The \textit{switching manifold} is defined as $\Sigma = \{ \vx \in \mR^n | \vv\T \vx = \delta \}$. 
\end{definition}
\begin{proposition}{(see \cite{Carmona2002})}
	All 2CPL$_n$ systems can be transformed to the \textit{Lure-like} form 
	\begin{equation}\label{eq_def_lurelike}
		\dot{\vx} = \vA \vx + \vb  + \vc \max(\ve_1\T \vx,0) = \begin{cases}
			\vA^- \vx + \vb, \quad \mbox{for } x_1 < 0 \\ 
			\vA^+ \vx + \vb, \quad \mbox{for } x_1 \geq 0
		\end{cases},
	\end{equation}
	where $\vA \in \mRnn, \vx,\vb,\vc \in \mR^n$, $\ve_1$ is the first basis vector of $\mR^n$ and the matrices $\vA^\pm$ being defined as $\vA^- = \vA $ and $\vA^+ = \vA + \vc \ve_1\T$. The switching manifold can thus be given as the hyperplane $x_1 = \ve_1\T \vx = 0$. When $\vb=0$, we speak of a \textit{homogeneous} 2CPL$_n$ system. The continuity property can be easily seen since $\vA^\pm$ only differ in their first column, and therefore, at $\Sigma = \{ \vx \in \mR^n | x_1=0\}$, both vector fields take the same value. \label{prop_lure}
\end{proposition}
\begin{remark}
	A general proof of Proposition \ref{prop_lure} is given in \cite{Carmona2002}. Here, we will prove the assertion for a class of driven piecewise linear mechanical systems. 
\end{remark}
\begin{remark}
	In \cite{Carmona2002}, 2CPL$_n$ systems are considered in the autonomous framework. Here, we use the term 2CPL$_n$ not only for autonomous systems but also for harmonically forced systems, which can be described by \eqref{eq_def_lurelike} with an additional harmonic forcing term.
\end{remark}
We consider a class of mechanical systems, which essentially consist of a linear multi-degree of freedom oscillator with symmetric system matrices, which is additionally equipped with a single unilateral spring and subjected to harmonic forcing. The unilateral spring models elastic contact within the system (e.g. a crack in a beam) or an elastic one-sided support. Such systems may be expressed by
\begin{align}\label{eq_eom_2ndorder}
	&\vM \ddot{\vq} + \alpha \vC \dot{\vq} + \vK \vq = \vw \lambda + \vf \cos(\Omega t),  \quad \mbox{bilinear multi-DOF oscillator}  \\ &\mbox{with } -\lambda = \max(k_n g, 0), \quad \mbox{and } g = \vw\T \vq - \delta \in \mR,\nonumber
\end{align}
where $\vq \in \mR^N$ is the generalized position vector, $\vw \in \mR^N$ is the vector of the generalized force direction of the contact force $\lambda$, $\vM=\vM\T \in \mRNN$ is the positive definite mass matrix, $\vC=\vC\T \in \mRNN$ is the damping matrix with $\alpha\geq0$ being an additional scalar parameter steering the damping in the system, $\vK=\vK\T \in \mRNN$ is the positive definite stiffness matrix, $\vf \in \mR^N$ is the external forcing vector and $g = \vw\T \vq - \delta \in \mR$ is the contact displacement. This equation describes a mechanical structure containing a unilateral spring which becomes active only if the contact displacement $g$ is positive. The nonlinearity has a nonsmooth nature despite being continuous, and is concentrated at the contact interface $g=0$. The characteristic of the unilateral elastic contact is shown in Figure \ref{fig_unilateral_elastic_force_law}.
The choice of a static contact force depending on the displacement is motivated by the study of the NNMs, but this choice is not a restriction for what follows and a velocity-dependent contact force can also be treated. 

In order to write the equations of motion \eqref{eq_eom_2ndorder} in first-order Lure-like form, we assume without loss of generality that the first component of the vector of generalized force directions $\vw$ is non-zero, as we are free to reorder the generalized positions. Therefore, it holds that $\vw = \begin{pmatrix}
	w_1 & \bar{\vw}\T 
\end{pmatrix}\T \in \mR^N$, with $w_1 \neq 0$ and $\bar{\vw} \in \mR^{N-1}$ arbitrary. We introduce the new set of coordinates $\tilde{\vq}$ and $\dot{\tilde{\vq}}$ defined as follows: 
\begin{equation}\label{eq_trafo_q_to_qtilde}
	\tilde{\vq} = \vV \vq  - \delta \tilde{\ve}_1, \quad \mbox{and} \quad \dot{\tilde{\vq}} = \vV \dot{\vq} ,
\end{equation}
with 
\begin{equation}
	\vV = \begin{pNiceArray}{cc}[ columns-width = 2cm]
		 w_1 & \bar{\vw}\T \\
		 \mathbf{0}_{(N-1) \times 1} & \vI_{(N-1)\times(N-1)} 
	\end{pNiceArray} , \nonumber
\end{equation}
where $\tilde{\ve}_1 \in \mR^N$ is the first basis vector of $\mR^N$ and $\vI$ denotes the identity matrix. Obviously, it holds that
\begin{equation}
	\tilde{\ve}_1\T \tilde{\vq} = g.
\end{equation}
Since $w_1 \neq 0$, the transformation matrix $\vV$ is invertible and the physical coordinates can be recovered from the new coordinates as follows:
\begin{equation}\label{eq_trafo_qtilde_to_q}
	\vq = \vV^{-1} (\tilde{\vq} + \delta \tilde{\ve}_1), \quad \mbox{and} \quad \dot{\vq} = \vV^{-1}\dot{\tilde{\vq}}. 
\end{equation}
Substituting these expressions in the equations of motion \eqref{eq_eom_2ndorder} and premultiplying by $\vV\negTranspose$ yields
\begin{equation}
	\tilde{\vM} \ddot{\tilde{\vq}} + \alpha \tilde{\vC} \dot{\tilde{\vq}} + \tilde{\vK}\tilde{\vq} + k_n \max(\tilde{\ve}_1\T \tilde{\vq} ,0) \tilde{\ve}_1 = \tilde{\vf} \delta + \tilde{\vf}_t  \cos(\Omega t)  
\end{equation}
with the transformed matrices and vectors defined as follows:
\begin{equation}\label{eq_transformed_matrices}
	\tilde{\vM} = \vV\negTranspose \vM \vV^{-1}, \quad \tilde{\vC} = \vV\negTranspose \vC \vV^{-1}, \quad \tilde{\vK} = \vV\negTranspose \vK \vV^{-1}, \quad \tilde{\vf} = - \vV\negTranspose \vK \vV^{-1} \tilde{\ve}_1, \quad \tilde{\vf}_t = \vV\negTranspose \vf. \nonumber
\end{equation}
By defining a state vector $\vx = \begin{pmatrix}
	\tilde{\vq}\T & \dot{\tilde{\vq}}\T 
\end{pmatrix}\T \in \mR^{n}$ with $n = 2N$, the system equations can be rewritten in first-order form as
\begin{equation}\label{eq_lurelike}
	\dot{\vx} = \vA(\alpha) \vx + \vc \max(\ve_1\T \vx, 0) + \vb + \vb_t \cos(\Omega t) = \begin{cases}
		\vA^-(\alpha) \vx + \vb \delta + \vb_t \cos(\Omega t), \quad \mbox{for } \ve_1\T \vx < 0\\
		\vA^+(\alpha) \vx + \vb \delta + \vb_t \cos(\Omega t), \quad \mbox{for } \ve_1\T \vx \geq 0\\
	\end{cases}, 
\end{equation}
with 
\NiceMatrixOptions{cell-space-limits = 5pt}
\begin{equation}\label{eq_A_c_b_bt_definition}
	\vA(\alpha) = \begin{pNiceArray}{cc}
		\bm{0}_{N\times N} & \vI_{N\times N} \\ -\tilde{\vM}^{-1} \tilde{\vK} &  -\alpha \tilde{\vM}^{-1} \tilde{\vC}
	\end{pNiceArray}, \quad \vc = \begin{pNiceArray}{cc}
	\bm{0}_{N\times 1} \\ -k_n \tilde{\vM}^{-1} \tilde{\ve}_1  
\end{pNiceArray} , \quad \vb = \begin{pNiceArray}{cc}
		\bm{0}_{N\times 1} \\ \tilde{\vM}^{-1} \tilde{\vf} 
	\end{pNiceArray}, \quad \vb_t = \begin{pNiceArray}{cc}
		\bm{0}_{N\times 1} \\  \tilde{\vM}^{-1} \tilde{\vf}_t  
	\end{pNiceArray}.
\end{equation}
The system matrices are given by $\vA^-(\alpha) = \vA(\alpha) $ and $\vA^+(\alpha) = \vA(\alpha) + \vc \ve_1\T$. 
The dynamics \eqref{eq_lurelike} is in the Lure-like form of \eqref{eq_def_lurelike} with an additional time-dependent forcing term. The system matrix $\vA$ is introduced as a function of the scalar $\alpha \geq 0$. In order to study the NNMs of the 2CPL$_n$ system, the underlying unforced and conservative dynamics has to be considered. This is easily achieved by setting $\alpha = 0$ and canceling the external forcing term in the original equations of motion \eqref{eq_lurelike}, which yields the following equation:
\begin{equation}\label{eq_conserv_auton_CPL}
	\dot{\vx} = \vA(0) \vx + \vc \max(\ve_1\T \vx, 0) + \vb \delta = \begin{cases}
		\vA_0^- \vx + \vb \delta, \quad \mbox{for } \ve_1\T \vx < 0\\
		\vA_0^+ \vx + \vb \delta, \quad \mbox{for } \ve_1\T \vx \geq 0\\
	\end{cases}. 
\end{equation}
with $\vA_0^- = \vA(0) $ and $\vA_0^+ = \vA(0) + \vc \ve_1\T$. 

\section{Invariant cones of homogeneous 2CPL$_n$ systems}
In this section, we will summarize the concept of invariant cones with a special focus on homogeneous 2CPL$_n$ systems, i.e. systems of the form \eqref{eq_def_lurelike} with $\vb = \bm{0}$. Originally, this concept was introduced by Küpper et al. \cite{Kuepper2008, Kuepper2011} for general, i.e. continuous and discontinuous, $n$-dimensional PWL systems. It is important to note that invariant cones have previously been observed for 2CPL$_3$ as invariant manifolds foliated by periodic orbits in the work by Carmona and coworkers \cite{carmona2005bifurcation, carmona2005limit, Carmona2006}. 
Consider homogeneous 2CPL$_n$ systems for which the dynamics reads as 
\begin{equation}\label{eq_homogeneousCPL}
	\dot{\vx} = \begin{cases}
		\vA^- \vx, \quad \mbox{for } x_1 < 0 \\ 
		\vA^+ \vx, \quad \mbox{for } x_1 \geq 0
	\end{cases},
\end{equation}
with the continuity condition on the switching hyperplane given by
\begin{equation}\label{eq_contcond_hcpl}
	\vA^+\vx^* = \vA^- \vx^*, \, \forall \vx^* \in \Sigma:= \{ \vx \in \mR^n \mid \ve_1\T \vx = x_1 = 0 \}.
\end{equation}
The state space can be decomposed in two subregions $\mathcal{A}^\pm$ separated by the switching hyperplane $\Sigma$ as follows
\begin{equation}
	\mR^n = \mathcal{A}^- \cup \Sigma \cup \mathcal{A}^+, \quad \mbox{with } \quad  \mathcal{A^\pm} := \{ \vx \in \mR^n \mid \ve_1\T \vx \gtrless 0 \}. \nonumber 
\end{equation}
In the following, we refer to the linear subsystems governing the dynamics in $\mathcal{A}^-$ and $\mathcal{A}^+$ by the $\ominus$ and $\oplus$ subsystem, respectively. As opposed to Filippov systems \cite{Lein2004}, systems of the form \eqref{eq_homogeneousCPL} with continuity condition \eqref{eq_contcond_hcpl} cannot exhibit sliding behavior on $\Sigma$, and thus, trajectories approaching $\Sigma$ cross it immediately. 
Let $\vph(t,t_0,x_0) := \vx(t)$ be the solution of \eqref{eq_homogeneousCPL} with initial condition $\vx(t_0)= \vx_0$. Depending on the direction of crossing, the switching hyperplane $\Sigma$ can be decomposed in the disjunct subsets 
\begin{equation}
	\Sigma^\pm := \{ \vx \in \Sigma \mid \exists \, \tau > 0 \, : \, \vph(t,0,\vx) \in \mathcal{A}^\pm, \, \forall \,  t \in (0,\tau) \},
\end{equation}
and the origin $\{\bm{0}\}$ being the equilibrium of the system, with $\Sigma = \Sigma^- \cup \Sigma^+ \cup \{\bm{0}\}$. The subsets $\Sigma^\pm$ are defined such that if $\vx \in \Sigma^\pm$, then the solution will continue in $\mathcal{A}^\pm$. Note that the subsets $\Sigma^\pm$ contain the half-hyperplanes 
\begin{equation}\label{eq_crossingmanifold}
	\begin{aligned}
		\{ \vx \in \Sigma \mid \ve_1\T \vA^\pm \vx < 0 \} \subset \Sigma^- \\
		\{ \vx \in \Sigma \mid \ve_1\T \vA^\pm \vx > 0 \} \subset \Sigma^+
	\end{aligned}
\end{equation}
but also part of the subset $\{ \vx \in \Sigma \mid \ve_1\T \vA^\pm \vx = 0 \}$ and are therefore cumbersome to express explicitly, e.g. 
\begin{align}
	\Sigma^- = &\{ \vx \in \Sigma \mid \ve_1\T \dot{\vx} < 0 \} \cup \{ \vx \in \Sigma \mid  \ve_1\T \dot{\vx} = 0, \,  \ve_1\T \ddot{\vx} < 0 \} \cup \cdots  \nonumber \\  &\cup  \{ \vx \in \Sigma \mid  \ve_1\T \vx^{(k)} = 0, \, \mbox{for } k = 1 \cdots (n-2), \,  \ve_1\T \vx^{(n-1)} < 0 \}.
\end{align}
We also note that both $\vA^-$ or $\vA^+$ can be used in \eqref{eq_crossingmanifold} without changing the definition, due to the continuity condition \eqref{eq_contcond_hcpl}.  
Furthermore, oscillatory mechanical systems are of interest in this work and therefore, we assume the presence of complex eigenvalues in both subsystems, which in turn introduces rotation in the state space and excludes trajectories confined to only one linear subsystem. Following this argumentation, a solution that evolves in one subregion and reaches $\Sigma$ necessarily crosses it towards the other subregion. The intersection point of the solution and the switching hyperplane can then be used as initial condition for the following part of the orbit. Hence, the solution of the homogeneous 2CPL$_n$ \eqref{eq_homogeneousCPL} can be described using exponential matrices as follows
\begin{equation}\label{eq_sol_composition}
	\vx(t) = \begin{cases}
		\exp(\vA^- (t - t_0)) \vx_0 , & t_0 \leq t \leq t_1 ,\\
		\exp(\vA^+ (t - t_1)) \exp(\vA^- (t_1 - t_0)) \vx_0, & t_1 \leq t \leq t_2 ,\\
		\cdots & \cdots
	\end{cases}
\end{equation}
where, without loss of generality, we assumed an initial point $\vx_0 = \vx(t_0) \in \Sigma^-$, such that a trajectory starting from $\vx_0$ leaves $\Sigma^-$ toward the region $\mathcal{A}^-$. Despite the apparent simplicity of \eqref{eq_sol_composition}, the difficulty lies in determining the crossing time instants $\{ t_1, t_2, \cdots \}$, which are \textit{a priori} unknown and are implicitly defined through $\vx(t_i)\in \Sigma, \, i = 1,2, \dots$ . The 2CPL$_n$ dynamics \eqref{eq_homogeneousCPL} can be studied by means of a discrete Poincaré map with a Poincaré section fixed at the switching hyperplane $\Sigma$. Assuming an initial condition ${\bm{\xi} = \vx(t_0) \in \Sigma^-}$, the orbit starting from $\bm{\xi}$ at $t_0$ reaches $\Sigma$ again after the time interval $t^- = t_1(\bm{\xi}) - t_0$ at the point $\bm{\eta} := \exp(\vA^- t^-) \bm{\xi} \in \Sigma^+$. The crossing occurs towards $\mathcal{A}^+$, and the orbit oscillates in the $\oplus$ subsystem until it reaches $\Sigma^-$ again after the time interval $t^+ = t_2(\bm{\eta}) - t_1(\bm{\xi})$ at the crossing point $\bm{\chi} := \exp(\vA^+ t^+) \bm{\eta} \in \Sigma^-$. In other words, we can define the Poincaré half-maps $\mathcal{P}^\pm$, whose combination yields the Poincaré map $\mathcal{P}$ of the whole system \eqref{eq_homogeneousCPL} as follows:
\begin{align}
	\mathcal{P}^- : \, \Sigma^- &\mapsto \Sigma^+ \nonumber \\
	\bm{\xi}  &\mapsto \bm{\eta} :=  \exp(\vA^- t^-) \bm{\xi}  \nonumber \\
	\mathcal{P}^+ : \, \Sigma^+ &\mapsto \Sigma^- \nonumber \\
	\bm{\eta}  &\mapsto \bm{\chi} :=  \exp(\vA^+ t^+ ) \bm{\eta} \nonumber \\
	\mathcal{P} := \mathcal{P}^+&(\mathcal{P}^-(\bm{\xi})) = \exp(t^+(\bm{\eta}) \vA^+) \exp(t^-(\bm{\xi})\vA^-) \bm{\xi},
\end{align}
where the first return time intervals $t^\pm$ can be defined as 
\begin{align}\label{eq_def_firstreturntimes}
	t^-(\bm{\xi}) &:= \inf \{ t> 0 \mid \ve_1\T \exp(t \vA^-) \bm{\xi}  = 0 , \quad \mbox{with } \bm{\xi} \in \Sigma^- \} , \nonumber \\
	t^+(\bm{\eta}) &:= \inf \{ t> 0 \mid \ve_1\T \exp(t \vA^+) \bm{\eta}  = 0, \quad \mbox{with } \bm{\eta} \in \Sigma^+  \} .
\end{align}
Due to the homogeneity property of \eqref{eq_homogeneousCPL}, solutions in both subregions can be given in closed form using the matrix exponential function. Scaling the initial conditions in \eqref{eq_def_firstreturntimes} does not change the definition of the first return time intervals, such that $t^\pm$ are positive homogeneous functions of degree zero, i.e. it holds that 
\begin{equation}
	t^-(\beta \bm{\xi}) = t^-(\bm{\xi}), \quad t^+(\beta \bm{\eta}) = t^+(\bm{\eta}), \quad \forall \beta > 0.
\end{equation}
Therefore, one can easily see that $t^\pm$ are constant on straight half-lines on $\Sigma$ emanating from the origin. This property extends also to the Poincaré half-maps $\mathcal{P}^\pm$ and the complete Poincaré map $\mathcal{P}$, which are homogeneous functions of degree one
\begin{equation}\label{eq_homog_prop_poincare}
	\mathcal{P}^- (\beta \bm{\xi}) = \beta \mathcal{P}^- ( \bm{\xi}), \quad \mathcal{P}^+ (\beta \bm{\eta}) = \beta \mathcal{P}^+ ( \bm{\eta}), \quad 
	\mathcal{P} (\beta \bm{\xi}) = \beta \mathcal{P} ( \bm{\xi}), \quad \forall \beta >0.
\end{equation}
\begin{figure}
	\centering	
	\includegraphics[scale=0.5]{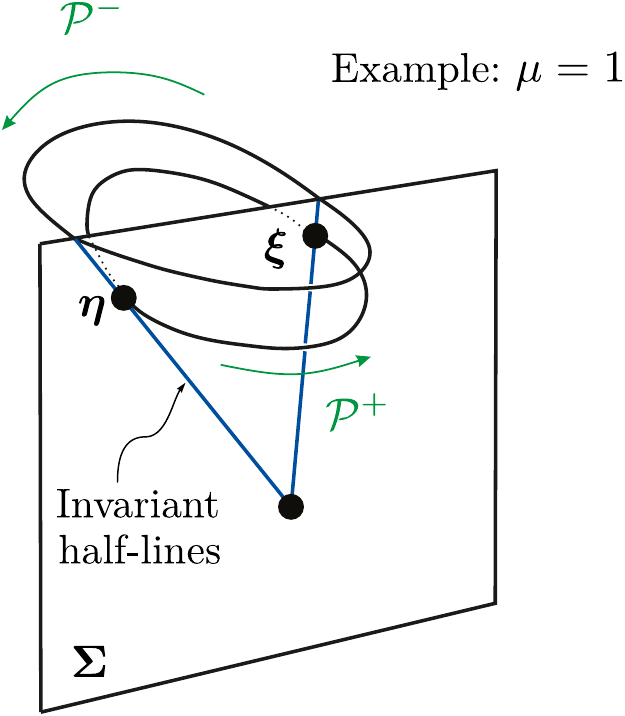}	
	\caption{A pictorial view of an invariant cone foliated by periodic orbits with $\mu=1$.}
	\label{fig_invcone_periodic}
\end{figure}%
\noindent%
For homogeneous 2CPL$_n$ systems of the form \eqref{eq_homogeneousCPL}, the equilibrium is located on the switching hyperplane and thus, the Jacobian is not unique. Therefore, the dynamics cannot be investigated using a linear analysis, since no "proper" eigenvalues are defined at the equilibrium. With the aim of reducing PWL systems to lower dimensional reduced systems containing their essential stability properties, Küpper \cite{Kuepper2008} introduced invariant cones as a nonsmooth counterpart to the center manifold for smooth systems. An invariant cone is defined as a special invariant set characterized by the existence of a straight half-line $\bm{\xi} \in \Sigma^-$ which fulfills the condition 
\begin{equation}\label{eq_def_invco_poincare}
	\mathcal{P}(\bm{\xi}) = \exp(t^+(\bm{\eta}(\bm{\xi})) \vA^+) \exp(t^-(\bm{\xi})\vA^-) \bm{\xi} = \mu \bm{\xi},
\end{equation}
where $\mu \in \mR^+$ is a positive real number. The magnitude of $\bm{\xi} \in \Sigma^-$ can be chosen as $\left\lVert \bm{\xi} \right\rVert = 1$ without any loss of generality. Similar arguments can be made if initial conditions are chose in $\Sigma^+$ by changing the order of the matrix exponential functions. 
Equation \eqref{eq_def_invco_poincare} can be interpreted as a \textit{nonlinear eigenvalue problem}, with a generalized eigenvalue $\mu$ and a generalized eigenvector $\bm{\xi}$. A solution of \eqref{eq_def_invco_poincare}, if it exists, is characterized by the four unknowns $(\bm{\xi},t^-, t^+, \mu) \in \mR^{n+3}$, which represent an invariant half-line, the two first return times and the generalized eigenvalue, respectively. This defines an invariant cone due to $\mathcal{P}(\beta \bm{\xi}) = \beta \mathcal{P}(\bm{\xi}) = \beta \mu \bm{\xi}, \, \forall \beta>0$.
This positive homogeneity property of the Poincaré map is essential for the following discussion and has an intuitive geometric interpretation. From \eqref{eq_homog_prop_poincare}, one can directly see that a point on an invariant half-line in $\Sigma$ that goes through the origin is mapped by $\mathcal{P}$ to a point on the same straight half-line. This is illustrated in Figure \ref{fig_invcone_periodic} for the case $\mu = 1$, where the invariant cone is foliated by periodic orbits. 
For an orbit starting from an initial condition $\bm{\xi}_0 \in \Sigma^-$ which belongs to the invariant cone, the time evolution can be tracked by applying the Poincaré map $\mathcal{P}$ recursively, such that the following holds:
\begin{align}
	\bm{\xi}_1 := \vx(t^- + t^+) &= \mu \bm{\xi}_0 , \nonumber \\ 
	\bm{\xi}_2 := \vx(2(t^- + t^+)) = \mu \bm{\xi}_1 &= \mu^2 \bm{\xi}_0, \nonumber\\ 
	&\cdots \nonumber\\
	\quad \bm{\xi}_k := \vx(k(t^- + t^+)) &= \mu^k \bm{\xi}_0.\nonumber
\end{align}
Thus, the stability of trajectories on the invariant cone is automatically determined by the eigenvalue~$\mu$. If~$\mu > 1$, trajectories belonging to the invariant cone spiral out of the origin, which represents the equilibrium of the system. Therefore, the equilibrium is unstable. If, however, $\mu<1$, then trajectories on the invariant cone spiral in towards the origin. In this case, stability of the origin is determined by the stability of the map $\mathcal{P}$ with respect to complimentary directions. Lastly, if $\mu = 1$, then the invariant cone is foliated by a family of periodic orbits characterized by the relationship $\mathcal{P}(\bm{\xi}) = \bm{\xi}$, and no statement can be given regarding the stability of the equilibrium at the origin. 
Note that the eigenvalue $\mu$ enables to describe the dynamics restricted to the invariant cone and does not address its attractivity. 
Differentiating the homogeneity condition \eqref{eq_homog_prop_poincare} yields
\begin{equation}
	\begin{aligned}
		\frac{\mathrm{d}}{\mathrm{d}\beta} (\mathcal{P}(\beta \vxi)) &\stackrel{!}{=} \frac{\mathrm{d}}{\mathrm{d}\beta} (\beta \mathcal{P}(\vxi)), \\ 
		\frac{\partial \mathcal{P}}{\partial \vxi}\Bigr|_{\beta \vxi} \vxi &\stackrel{!}{=} \mathcal{P}(\vxi). 
	\end{aligned}
\end{equation} 
Therefore, if $\overline{\bm{\xi}}$ is a solution of the generalized eigenvalue problem \eqref{eq_def_invco_poincare} with eigenvalue $\mu$, then due to the relationship 
\begin{equation}
	\mathcal{P} (\overline{\bm{\xi}}) = \frac{\partial \mathcal{P}}{\partial \bm{\xi} }(\overline{\bm{\xi}}) \overline{\bm{\xi}} = \mu \overline{\bm{\xi}},
\end{equation}
one can see that $\overline{\bm{\xi}}$ is also an eigenvector of the corresponding Jacobian of the Poincaré map. If we assume that the remaining $(n-2)$ eigenvalues $\lambda_1, \cdots, \lambda_{n-2}$ of the Jacobian $\frac{\partial \mathcal{P}}{\partial \bm{\xi} }(\overline{\bm{\xi}})$ satisfy 
\begin{equation}\label{eq_attractivity}
	|\lambda_j| < \min(1,\mu), \quad j = 1, \cdots, n-2 ,
\end{equation}
then the invariant cone is attractive under the flow of \eqref{eq_homogeneousCPL} and the dynamics on the cone is governed by~$\mu$. In case the cone is foliated by trajectories spiraling in towards the origin, these trajectories converge faster to the cone than to the origin \cite{Weiss2012}.
In order to obtain the invariant cone numerically, a system of nonlinear equations has to be set up. As discussed previously, the central condition for the existence of an invariant half-line on the switching hyperplane determining the invariant cone is \eqref{eq_def_invco_poincare}, which involves the quantities ${\bm{\xi},\mu, t^-}$ and~$t^+$. Furthermore, additional conditions have to be incorporated to fix the eigenvector and its first return point to the switching hyperplane.
The nonlinear system to determine the invariant cone of the homogeneous 2CPL$_n$ system \eqref{eq_homogeneousCPL} reads as follows:
\begin{equation}\label{eq_classic_inv_cone_problem}
	\vF(\bm{\xi},t^-,t^+, \mu) = \begin{pmatrix}
		e^{t^+ \vA^+} e^{t^- \vA^-} \bm{\xi} - \mu \bm{\xi} \\
		\ve_1\T e^{t^- \vA^- } \bm{\xi} \\ 
		\ve_1\T \bm{\xi} \\ 
		\bm{\xi}\T \bm{\xi} - 1 
	\end{pmatrix}\stackrel{!}{=} \bm{0} .
\end{equation} 
Herein, the first equation is the explicit expression of the Poincaré map given in \eqref{eq_def_invco_poincare}, the second and third equations fix the first crossing point and the eigenvector to $\Sigma$, respectively, and the fourth equation is a normalization that can be chosen arbitrarily. This nonlinear system contains ${(n+3)}$ equations for ${(n+3)}$ unknowns and can be solved using a Newton scheme. Furthermore, an arclength-continuation algorithm can be used to trace families of invariant cones depending on a system parameter. 
 
\section{NNM backbones as extended invariant cones of inhomogeneous 2CPL$_n$ systems}\label{section:backbones}
In this section, we propose an augmented invariant cone problem to compute NNMs of inhomogeneous 2CPL$_n$ systems. Herein, the NNMs are represented by frequency-energy plots, also called "backbone" curves. The definition of NNMs used in this work is the extended version of Rosenberg´s concept, i.e. NNMs are understood as not necessarily synchronous periodic motions of the undamped autonomous system \cite{Peeters2009}. In this work, we determine only fundamental NNMs, i.e. those that are continuations of classical LNMs. Furthermore, the NNMs that can be computed by our modified invariant cone approach are characterized by one crossing of the switching hyperplane during one response period, which is inherent to the expression of the Poincaré map \eqref{eq_def_invco_poincare}. We hereby note that NNMs, in general, may bifurcate leading to the birth of so-called \textit{superabundant} NNMs without linear counterparts. Their modal curves in configuration space  may have more complex shapes implying periodic motions with a varying number of crossings per period. For instance, such NNMs have been observed in PWL systems such as in \cite{Chen1996, Chati1997, Jiang2004, Vestroni2008, Casini2011, Casini2012}. In fact, the possibility of NNMs becoming unstable through bifurcation mechanisms and the appearance of new more complex superabundant NNMs represents one of the fundamental differences between NNMs and the classical linear modes. In this work, the focus is placed on tracing fundamental NNM branches and determining their stability. The bifurcating superabundant modes that can emerge due to period doubling bifurcations or internal resonances are not considered here. 
Another key characteristic of nonlinear systems is the energy-dependence of resonance frequencies, a significant feature from an engineering perspective and generally present in all smooth nonlinear systems. However, in 2CPL$_n$ systems, this energy-dependence is observed only in inhomogeneous systems. Homogeneous systems, where the nonlinearity is fixed to the switching hyperplane, exhibit energy-independent NNMs similar to linear systems. This behavior has been observed in several studies, including \cite{Zuo1994}, \cite{Vestroni2008, Casini2012} and more recently in \cite{Attar2017}. For homogeneous 2CPL$_n$ systems, energy-independent NNMs can be regarded as periodic orbits foliating the corresponding invariant cones, if they cross the switching hyperplane once per period.
However, in the presence of an inhomogeneous term, i.e. $\vb \neq \bm{0}$ in the Lure-like form \eqref{eq_def_lurelike}, the NNMs become energy dependent, which is obvious from a physical perspective. In fact, due to the presence of a nonvanishing clearance gap between the structure and the unilateral contact, the motion may remain purely linear if the energy is not sufficient to make the structure interact with the constraint. In this scenario, the 2CPL$_n$ oscillates only in one subregion, displaying linear behavior. At higher energy levels, the oscillation amplitudes increase, causing solutions to leave the domain of one linear subsystem and cross the switching hyperplane. This introduces a stiffness change into the motion, responsible for shifting the overall frequency of the composed motion. This energy-dependence due to inhomogeneous terms contrasts with the formulation of invariant cones, which rely on positive homogeneity properties. In this section, we propose a modified invariant cone formulation that accounts for the inhomogeneity by introducing an augmented homogeneous higher-dimensional system. This allows the computation of classical invariant cones, from which specific solutions are interpreted as part of energy-dependent NNMs of the original inhomogeneous system dynamics. This modified invariant cone approach is then shown to yield a special form of the shooting method, which can be combined with a pseudo-arclength continuation procedure. The shooting method with path-following is outlined in Appendix A in a slightly modified form of \cite{Kerschen2014}, similar to the approach applied in \cite{Attar2017}. 
\subsection{Augmented invariant cone}
In the following, we investigate the unforced dynamics of \eqref{eq_eom_2ndorder} by canceling the external forcing term $\vf=\bm{0}$, which can be cast into first order Lure-like form \eqref{eq_lurelike} with $\vb_t=0$ yielding
\begin{equation}\label{eq_lurelike_autonomoussys}
		\dot{\vx} = \begin{cases}
		\vA^-(\alpha) \vx + \vb \delta, \quad \mbox{for } \ve_1\T \vx < 0\\
		\vA^+(\alpha) \vx + \vb \delta, \quad \mbox{for } \ve_1\T \vx \geq 0\\
	\end{cases}.
\end{equation}
We hereby note that the conservative dynamics are recovered for $\alpha = 0$ but we keep $\alpha$ as a parameter for now. The resulting system is an autonomous inhomogeneous 2CPL$_n$ system. In order to apply the invariant cone concept, the state is augmented by taking the gap length $\delta$ as an additional state denoted by $v_\delta$ with the aim of transforming the system to a homogeneous 2CPL$_{n+1}$ system. The dynamics can thus be rewritten as homogeneous system
\begin{equation}\label{eq_augmented_hom_sys_autonomous}
	\dot{\vv} = \tilde{\vA}^\pm (\alpha) \vv = \begin{cases}
	\tilde{\vA}^-(\alpha) \vv, \quad \mbox{for } \hat{\ve}_1\T \vv < 0\\
	\tilde{\vA}^+(\alpha) \vv, \quad \mbox{for } \hat{\ve}_1\T \vv \geq 0\\
	\end{cases},
\end{equation}
where the new augmented state vector is defined as ${\vv\T = \begin{pmatrix}
	\vx\T & v_\delta 
	\end{pmatrix} \in \mR^{n+1}}$ and $\hat{\ve}_1$ is the first basis vector of $\mR^{n+1}$. The augmented system matrices read as
\begin{equation}\label{eq_definition_Atilde}
	\tilde{\vA}^\pm (\alpha) = \begin{pNiceArray}{cc}
		\vA^\pm(\alpha) & \vb \\ \bm{0}\T & 0
	\end{pNiceArray} \in \mR^{(n+1)\times (n+1)},
\end{equation}
where $\vb$ is defined in \eqref{eq_A_c_b_bt_definition}. Note that due to the zeros in the last row of $\tilde{\vA}^\pm$, the additional state $v_\delta$ has the dynamics $\dot{v}_\delta= 0$ and thus remains constant. Hence, in order to recover the dynamics of system \eqref{eq_lurelike_autonomoussys}, $v_\delta$ must be initialized such that $v_\delta(t_0)=\delta$. For the homogeneous and autonomous augmented system \eqref{eq_augmented_hom_sys_autonomous}, the state space can be decomposed in two subregions $\mathcal{A}^\pm$ separated by the switching hyperplane $\Sigma$ as follows
\begin{equation}
	\mR^{n+1} = \mathcal{A}^- \cup \Sigma \cup \mathcal{A}^+, \quad \mbox{with } \quad  \mathcal{A^\pm} := \{ \vv \in \mR^{n+1} | \hat{\ve}_1\T \vv \gtrless 0 \}, \quad \mbox{and } \quad \Sigma := \{ \vv \in \mR^{n+1} | \hat{\ve}_1\T \vv = q_1 = 0 \}. \nonumber 
\end{equation} 
Let $\vph(t,t_0,\vv_0) := \vv(t) $ be the solution of \eqref{eq_augmented_hom_sys_autonomous} with initial condition $\vv(t_0)=\vv_0$. We assume furthermore that initial conditions are chosen as $ \vv_0 \in \Sigma^- $ with $\Sigma^-$ a subset of $\Sigma$ defined as
\begin{equation}
	\Sigma^- := \{ \vv \in \Sigma \mid \exists \, \tau > 0 \, : \, \vph(t,0,\vv) \in \mathcal{A}^-, \, \forall \,  t \in (0,\tau) \}.
\end{equation}
As discussed in Section $3$, an invariant cone problem can be formulated as in \eqref{eq_classic_inv_cone_problem}. Since we are interested in the conservative dynamics, all motions on the individual NNM manifolds must be periodic and $\alpha$ must fulfill $\alpha=0$. Thus, we seek invariant cones foliated by periodic orbits by setting $\mu = 1$ in \eqref{eq_classic_inv_cone_problem}. These invariant cones are obtained as half-lines on the switching hyperplane $\Sigma$, which represent a set of initial conditions ensuring periodic response. Moreover, the conservative dynamics implies that the total mechanical energy remains constant over the time period of the oscillation and its value is only determined by the initial condition. Since solutions of \eqref{eq_classic_inv_cone_problem} are a set of initial conditions located on the switching hyperplane $\Sigma$ and ensuring periodicity, the total energy of the periodic motion can be described by its expression in the $\ominus$ subsystem, which is much simpler than in the $\oplus$ subsystem. Therefore, the total energy for points on the switching hyperplane is obtained in the original coordinates $\vy\T = \begin{pmatrix}
\vq\T & \dot{\vq}\T 
\end{pmatrix}$ as 
\begin{equation}\label{eq_energy_minus_q}
E(\vy) = \frac{1}{2} \vy\T \begin{pmatrix}
\vK & \bm{0} \\ \bm{0} & \vM 
\end{pmatrix} \vy = a, \quad \dot{E}(\vy) = 0, \quad \mbox{for } \alpha = 0, 
\end{equation}
where $a$ denotes the constant energy level. Since the invariant cone problem is formulated in terms of the transformed coordinates $\vx\T = \begin{pmatrix}
\tilde{\vq}\T & \dot{\tilde{\vq}}\T
\end{pmatrix}$, the transformation \eqref{eq_trafo_qtilde_to_q} can be substituted in \eqref{eq_energy_minus_q} to yield the expression of the total energy in the new coordinates as follows:
\begin{equation}\label{eq_energyexpression}
	\begin{aligned}
		\tilde{E}^- (\tilde{\vq},\dot{\tilde{\vq}},\delta) &= \frac{1}{2} \left( \vV^{-1} \left( \tilde{\vq} + \delta \ve_1\right) \right)\T \vK \left( \vV^{-1} \left( \tilde{\vq} + \delta \ve_1\right) \right) + \frac{1}{2} \left( \vV^{-1} \dot{\tilde{\vq}} \right)\T \vM \left( \vV^{-1} \dot{\tilde{\vq}} \right) \\ 
		&= \frac{1}{2} \left( \tilde{\vq} + \delta \ve_1\right)\T \tilde{\vK} \left( \tilde{\vq} + \delta \ve_1\right) + \frac{1}{2} \dot{\tilde{\vq}}\T \tilde{\vM} \dot{\tilde{\vq}}.
	\end{aligned}
\end{equation} 
The resulting nonlinear problem determining the invariant cones foliated by periodic orbits for the homogeneous system \eqref{eq_augmented_hom_sys_autonomous} is then obtained as
\begin{equation}\label{eq_augmented_invcone_autonomous}
		\vF(\vv,t^-,t^+, \alpha, a) = \begin{pmatrix}
		e^{t^+ \tilde{\vA}^+(\alpha)} e^{t^- \tilde{\vA}^-(\alpha)} \vv - \vv \\
		\hat{\ve}_1\T e^{t^- \tilde{\vA}^-(\alpha) } \vv \\ 
		\hat{\ve}_1\T \vv \\ 
		\tilde{E}(\vv) - a  
	\end{pmatrix}\stackrel{!}{=} \bm{0} ,
\end{equation}
where the arbitrary normalization in the last row of \eqref{eq_classic_inv_cone_problem} is replaced by an energy constraint used to fix the value of the energy $\tilde{E}(\vv)$ to an unknown level $a$, with $\tilde{E}$ expressed as function of the augmented state $\vv$. 
The first equation in \eqref{eq_augmented_invcone_autonomous} is the explicit expression of the Poincaré map of the augmented system \eqref{eq_augmented_hom_sys_autonomous}, denoted here by $\mathcal{P}_a$. Following similar arguments as discussed in the previous section, the homogeneity property holds for $\mathcal{P}_a$ and can be stated for a vector $\vv\T = \begin{pmatrix}
	\vx\T & v_\delta
\end{pmatrix} \in \Sigma^-$ as follows: 
\begin{equation}
	\mathcal{P}_a(\beta \vv ) = \beta e^{t^- \tilde{\vA}^-(\alpha)} e^{t^+ \tilde{\vA}^+(\alpha)} \vv = \beta \mathcal{P}_a (\vv) = \beta \begin{pmatrix}
		\vx \\ v_\delta
	\end{pmatrix}, \quad \forall \beta > 0.
\end{equation}
The invariant cone problem \eqref{eq_augmented_invcone_autonomous} is therefore represented by a system of $(n+4)$ equations with $(n+5)$ unknowns $(\vv,t^-,t^+, \alpha, a)$. A solution of \eqref{eq_augmented_invcone_autonomous} denoted by $(\vv_s,t_s^-,t_s^+, \alpha_s, a_s)$  represents an invariant cone characterized by an invariant half-line on $\Sigma^-$ emanating from the origin and passing through the point $\vv_s \in \Sigma^-$. This invariant half-line can thus be described by $\beta \vv_s \in \Sigma^-, \, \forall \beta>0$. The invariant cone in the augmented state space $\mR^{n+1}$ is foliated by periodic orbits starting on the invariant half-line, all sharing a common constant period time obtained as $T_s = t_s^- + t_s^+$.  Moreover, the periodicity also implies conservative dynamics, and thus solutions must satisfy $\alpha_s=0$ in order to cancel the damping terms in $\tilde{\vA}^\pm$. The energy level $a_s$ corresponds to only one of all periodic orbits on the cone. In fact, the energy constraint in \eqref{eq_augmented_invcone_autonomous} is used as a constraint equation with the aim of obtaining a unique solution of the problem. Thus, all periodic orbits share the same period time $T_s$ regardless of their energy level. Here, it is important to highlight that the cone property described by the positive scaling factor $\beta>0$ applies to the augmented state vector $\vv_s$, thereby also scaling its last component $v_{\delta,s} = \hat{\ve}_{n+1}\T \vv_s$. This can be interpreted as an augmented invariant cone corresponding to the augmented 2CPL$_{n+1}$ system \eqref{eq_augmented_hom_sys_autonomous}. This augmented cone is, in turn, foliated by periodic solutions of a $1-$parameter family of 2CPL$_n$ systems of the form \eqref{eq_lurelike_autonomoussys}, which are parametrized by different clearance gap lengths given by $\beta v_{\delta,s} > 0, \, \forall \, \beta >0$, but all sharing the same free response period time. From all the periodic orbits forming the augmented invariant cone, only the one corresponding to the actual gap length $\delta$ of the original system \eqref{eq_lurelike_autonomoussys} is of interest. Therefore, from all possible values of the scaling parameter $\beta$, only one specific value $\beta^*$ can be chosen such that the periodic motion with the correct gap length $\delta$ of the original system is recovered:
\begin{equation}\label{eq_constraining_gaplength}
	\beta^* v_{\delta,s} \stackrel{!}{=} \delta, \quad \implies \beta^* = \frac{\delta}{v_{\delta,s}} > 0.
\end{equation}
Using $\beta^*$ as a scaling factor, the initial condition on the switching hyperplane characterizing the periodic motion of the original system \eqref{eq_lurelike_autonomoussys} can be obtained as 
\begin{equation}
	\hat{\vv}_p = \beta^* \vv_s,
\end{equation}
where $\vv_s$ is part of the solution of \eqref{eq_augmented_invcone_autonomous}. 
To summarize, solutions of \eqref{eq_augmented_invcone_autonomous} are invariant cones characterized by invariant half-lines of the Poincaré map $\mathcal{P}_a$ belonging to the augmented system \eqref{eq_augmented_hom_sys_autonomous}. These invariant cones are foliated by periodic orbits from which only one solution is of interest, which can be obtained by scaling $\vv_s$ using the factor $\beta^*$. Thus, this can be interpreted as a solution of \eqref{eq_augmented_invcone_autonomous} with an additional constraint equation. The resulting problem can be written in the form
\begin{equation}\label{eq_constrained_invcone}
	\tilde{\vF}(\vv,t^-,t^+, \alpha, a) = \begin{cases}
		\vF(\vv,t^-,t^+, \alpha, a) = \bm{0}, \\
		\hat{\ve}_{n+1}\T \vv = \delta ,
	\end{cases}
\end{equation}
where the first equation is the unconstrained invariant cone problem of the augmented system and the second equation represents the additional constraint to recover the correct dynamics with the gap length $\delta$.
In order to obtain a simplified form of \eqref{eq_constrained_invcone}, the constraint equation can be substituted into the unconstrained invariant cone problem. 
The second and third equations of \eqref{eq_augmented_invcone_autonomous} are not affected by the constraint $v_\delta = \delta$ since they only include the first component of $\vv$ which is the same as the first component of $\vx$. 
Substituting $\delta$ as a parameter into the expression of the total energy and considering $\vx$ as the only variable yields $\tilde{E}(\vx_0,\delta) - a$ as the fourth equation of \eqref{eq_constrained_invcone}. Finally, the constraint equation $v_\delta = \delta$ is substituted into the explicit expression of the Poincaré map. Due to the special form of the matrices $\tilde{\vA}^\pm(\alpha)$ defined in \eqref{eq_definition_Atilde}, we make use of the following special relation 
\begin{equation}\label{eq_saad_equation}
	\exp\left(\begin{bNiceArray}{cc}
		\vB & \vc \\ 0 & 0
	\end{bNiceArray}\right) = \begin{bNiceArray}{cc}
	e^{\vB} & \varphi_1(\vB)\vc \\ 0 & 1 
\end{bNiceArray},
\end{equation}
which was derived by Saad in \cite{Saad1992}, for a matrix $\vB \in \mathbb{C}^{n\times n}$, a vector $\vc \in \mathbb{C}^n$ and the \textit{phi-}function defined such that 
\begin{equation}
	\vB \varphi_1(\vB) = e^\vB - \bm{I}.
\end{equation}
Applying this relation to the matrices $\tilde{\vA}^\pm(\alpha)$ defined in \eqref{eq_definition_Atilde} yields the following equation
\begin{equation}\label{eq_poincaremap_augmented_withconstraint}
	e^{t^- \tilde{\vA}^-} e^{t^+ \tilde{\vA}^+} \vv - \vv = \begin{pmatrix}
		e^{t^+ \vA^+}\left( e^{t^- \vA^-} \vx + \varphi_1(t^-\vA^-) \tilde{\vb} \delta \right) +  \varphi_1(t^+\vA^+) \tilde{\vb} \delta - \vx \\ 0
	\end{pmatrix},
\end{equation}
where the dependency of $\tilde{\vA}^\pm$ and $\vA^\pm$ on $\alpha$ was omitted to shorten the notation. 
Obviously, the last row is always $0$ and can thus be discarded in \eqref{eq_augmented_invcone_autonomous}. 
By reforming the upper part of \eqref{eq_poincaremap_augmented_withconstraint}, we end up with the simplified equations of the constrained invariant cone problem which reads as:
\begin{equation}\label{eq_modified_invcone_autonomous}
		\vF(\vx_0, t^-, t^+, \alpha, a) = \begin{bmatrix}
			\begin{bmatrix}
				\vI_{n} & 0
			\end{bmatrix} e^{t^+ \tilde{\vA}^+(\alpha)}  e^{t^- \tilde{\vA}^-(\alpha)} \begin{bmatrix}
					\vx_0  \\ \delta
				\end{bmatrix} - \vx_0 \\
			\ve_1\T \vx_0 \\ 
			\ve_1\T \begin{bmatrix}
				\vI_{n} & 0
			\end{bmatrix} e^{t^- \tilde{\vA}^-(\alpha)} \begin{bmatrix}
				\vx_0 \\ \delta
			\end{bmatrix} \\ 
			\tilde{E}(\vx_0,\delta) - a 
		\end{bmatrix}.
\end{equation}
This nonlinear problem contains $(n+3)$ equations for $(n+4)$ unknowns, which in turn is a classic continuation problem by regarding the energy level $a$ as the continuation parameter. Solutions of \eqref{eq_modified_invcone_autonomous} trace a branch in the $(n+4)$-dimensional space of initial conditions $\vx_0$, return times $t^\pm$, the damping parameter $\alpha$ and the energy levels $a$. Due to the structure of the matrices $\vA^\pm(\alpha)$, periodic motions are only possible if $\alpha=0$, which has to hold along the whole branch of periodic orbits. 
By applying a continuation algorithm on problem \eqref{eq_modified_invcone_autonomous}, the NNM backbone curves can be traced. This is equivalent to solving the augmented invariant cone problem \eqref{eq_augmented_invcone_autonomous} at each energy level, and then constraining the solution to recover the correct periodic orbit corresponding to the gap length $\delta$ of the original system  \eqref{eq_lurelike_autonomoussys}. This is repeated along the branch at different energy levels, thereby recovering the frequency-energy dependence of the NNMs. Furthermore, we should highlight that solutions of \eqref{eq_modified_invcone_autonomous} must start on the switching hyperplane and exhibit crossing behavior such that linear oscillations restrained to one subsystem cannot be computed. This implies that \eqref{eq_modified_invcone_autonomous} describes only the nonlinear regions of the backbone curves, for which frequencies depend on the energy level, whereas the energy-independent parts of the backbones characterized by straight lines at constant frequency cannot be traced. This, however, is not problematic, since these parts of the backbones correspond to the classic LNMs and are fixed at the underlying linear frequencies of the $\ominus$ subsystem, which can be easily computed. For each NNM, the frequency-energy backbone curve is a straight line up to the point where the corresponding periodic solution crosses the switching hyperplane $\Sigma$. After that point, the nonsmooth nonlinearity is activated, thereby introducing the frequency-energy dependence in the backbone curve. It is also important to point out that the structure of \eqref{eq_modified_invcone_autonomous} is reminiscent of the shooting method presented in Appendix \ref{Appendix_shooting}. In fact, the first equation in \eqref{eq_modified_invcone_autonomous} represents a periodicity condition obtained by the time integration of the system over one period and comparing it to the initial condition. The oscillation period is implicitly assumed to include two intervals $t^\pm$ corresponding to oscillations in both subsystems. The constraint equations are chosen such that solutions are guaranteed to cross the switching hyperplane, thereby excluding any solutions confined to the $\ominus$ subsystem. Therefore, we can establish that the invariant cone problem of the augmented system constrained to the correct gap length of the inhomogeneous system reduces to a special form of the shooting method. Solutions of this modified invariant cone problem trace the nonlinear regions of the NNM backbone curves, which are characterized by a crossing behavior. 
\subsection{Stability}\label{Subsection_Stability}
The stability properties of NNMs of 2CPL$_n$ system have been investigated in a large amount of papers dealing with both homogeneous and inhomogeneous systems, such as \cite{Zuo1994, Jiang2004, Casini2012} and \cite{Attar2017}. In the case of conservative mechanical systems, stability of periodic orbits is understood in the sense of Lyapunov and is studied by means of two methods: Floquet multipliers and Poincaré maps. The Floquet multipliers are obtained as the eigenvalues of the monodromy matrix, which is governed by the homogeneous part of the 2CPL$_n$ system and an initial condition given by the identity matrix:
\begin{equation}
	\dot{\vPh} = \vA^\pm \vPh, \quad \vPh(t=0) = \vI_{n\times n}. 
\end{equation}
The stability can be inferred from the magnitudes of the eigenvalues of $\vPh$. A necessary condition for stability is that all the Floquet multipliers lie on or inside the unit circle in the complex plane. Thus, the Floquet multipliers can be used to detect bifurcation points, at which branches of NNMs exhibit stability changes. As observed in \cite{Jiang2004} in the context of an inhomogeneous 2CPL$_2$ system, the fundamental NNMs are stable up to a bifurcation point where two Floquet multipliers merge at $-1$ with one of them leaving the unit circle. Beyond this critical point, quasi-periodic and chaotic behavior have been observed. The fundamental NNM branches were also observed to regain stability at another critical point, restoring the stability of the periodic behavior. In general, quasi-peridic behavior beyond the instability of fundamental branches is studied by means of Poincaré maps \cite{Zuo1994, Jiang2004, Casini2012, Attar2017}. The corresponding Poincaré sections have to be chosen appropriately in order to capture the complex behavior. These Poincaré maps can be constructed numerically given that solutions are linear in each subregion.
In this work, our focus lies on determining fundamental NNM branches using the modified invariant cone formulation \eqref{eq_modified_invcone_autonomous}, and post-critical responses, i.e. bifurcating NNMs and/or quasi-periodic motions, are not considered here. The stability of the NNMs obtained by our formulation follows as a by-product of the solution of the nonlinear problem \eqref{eq_modified_invcone_autonomous}. In fact, since the system is continuous, the saltation matrices needed to compute the monodromy matrix of the general 2CPL$_n$ system are the identity matrices \cite{Lein2004}. Therefore, the monodromy matrix can be obtained automatically using the matrix exponential function along with the first return times $t^\pm$ and $\alpha$ obtained from the solution of \eqref{eq_modified_invcone_autonomous} as follows:
\begin{equation}\label{eq_monodromy_expm}
	\vPh = e^{t^+ \tilde{\vA}(\alpha)^+} e^{t^- \tilde{\vA}(\alpha)^-}.
\end{equation}
Therefore, upon solving \eqref{eq_modified_invcone_autonomous} using a path-following algorithm, stability can be evaluated at each point of the branch by inspecting the magnitudes of the eigenvalues of \eqref{eq_monodromy_expm}. 
\subsection{Numerical example}
\begin{figure}
	\centering	
	\includegraphics{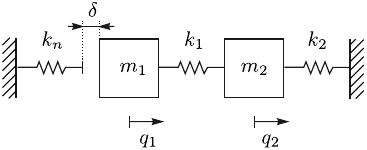}	
	\caption{A 2CPL$_2$ mechanical system with a unilateral elastic contact.}
	\label{fig_mech_sys_jiang2004}
\end{figure}

\begin{figure}[t]
	\begin{subfigure}{0.49\linewidth}
		\centering
		\includegraphics{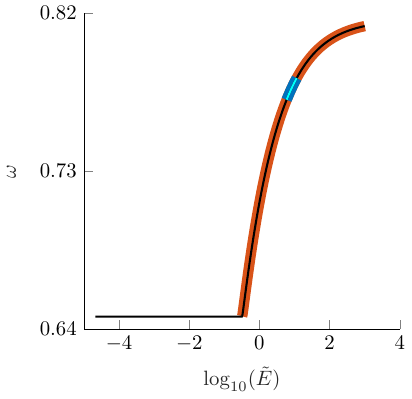}%
		\caption{}
		\label{fig:backbone_mode1}
	\end{subfigure}
	\begin{subfigure}{0.49\linewidth}
		\centering
		\includegraphics{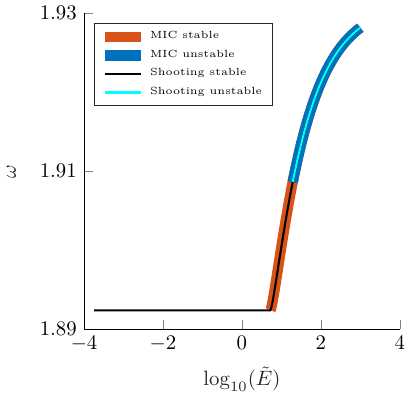}
		\caption{}
		\label{fig:backbone_mode2}
	\end{subfigure}
	\caption{Backbone curves showing the validity of the MIC approach by comparison to classical shooting.}
	\label{fig:backbones}
\end{figure}

We investigate the 2CPL$_4$ system studied by Jiang et al. \cite{Jiang2004} to illustrate how our modified invariant cone (MIC) formulation captures the nonlinear backbone curves and the NNM invariant manifolds. The system is sketched in Figure \ref{fig_mech_sys_jiang2004}. 
In order to satisfy the assumption of a nonzero first component in the vector of generalized force directions $\vw$ stated in Section 2, we reorder the set of generalized positions chosen in \cite{Jiang2004} and use the set of generalized coordinates $\vq =  \begin{pmatrix}
q_1 & q_2 
\end{pmatrix}\T $ as shown in Figure \ref{fig_mech_sys_jiang2004} to express the system dynamics. 
The equations of motion of the system in these coordinates are 
\begin{equation}\label{eq_Jiang2004_eom}
	\begin{pmatrix}
	m_1 & 0 \\ 0 & m_2 
	\end{pmatrix} \ddot{\vq} + \begin{pmatrix}
	k_1 & -k_1 \\ -k_1 & k_1+k_2 
	\end{pmatrix} \vq = \vw \lambda , \quad \mbox{with } -\lambda = \max(k_n g , 0), \quad \mbox{and } g = \vw\T \vq - \delta,
\end{equation}
where the generalized force direction is given by $\vw = \begin{pmatrix}
	-1 & 0
\end{pmatrix}\T $ and $\delta$ is the clearance gap length. 
In order to apply the invariant cone approach, the new set of coordinates $\tilde{\vq}$ is introduced as defined in \eqref{eq_trafo_q_to_qtilde}. The system can be cast in first order Lure-like form \eqref{eq_lurelike_autonomoussys} using the new state vector $\vx\T = \begin{pmatrix}
\tilde{\vq}\T & \dot{\tilde{\vq}}\T
\end{pmatrix}$, where the system matrices and the inhomogeneous term are defined in \eqref{eq_transformed_matrices} and \eqref{eq_A_c_b_bt_definition}. Note that the system \eqref{eq_Jiang2004_eom} does not include damping. Nevertheless, the parameter $\alpha$ is important for the invariant cone problem, and therefore, we set the matrix $\tilde{\vC}=\tilde{\vM}$ in \eqref{eq_A_c_b_bt_definition}, knowing that the solution of the modified invariant cone problem \eqref{eq_modified_invcone_autonomous} must satisfy $\alpha=0$ to obtain the conservative dynamics investigated here. 
We will compute the NNMs of the system with the parameter values used by Jiang et al. \cite{Jiang2004},
\begin{equation}\label{eq_Jiang_paramvalues}
	m_1=m_2=1, \quad k_1 = 1.5,\quad k_2=1 ,\quad k_n = 1.5, \quad \delta = 1.
\end{equation} 
The linearized eigenfrequencies of the system without unilateral spring are ${\omega_1=0.6472}$ and ${\omega_2=1.8924}$. 
For the NNMs emanating from the linear modes of the underlying linear dynamics, the backbone curves are traced by solving the nonlinear problem \eqref{eq_modified_invcone_autonomous} combined with a continuation algorithm. An initial point on the branch can be easily obtained by computing the linear mode and finding the point at which solutions start crossing the switching hyperplane $\Sigma$, and then following the branch in the direction of increasing energy. Alternatively, one can start at a high energy level and find an initial point by a Newton iteration; from there the backbone is traced backwards in the direction of decreasing energy. To compare the results, the classic shooting method with path-following is used, which is summarized in Appendix \ref{Appendix_shooting}. The nonlinear frequency is obtained from the first return times $t^\pm$ by $\omega= 2\pi / (t^-+t^-)$. The resulting NNM backbone curves are shown in Figure \ref{fig:backbones}. For low energy levels, solutions oscillate in the LNM without crossing the switching hyperplane and the backbone curves are straight lines fixed at the two linear eigenfrequencies. These energy-independent linear modes cannot be obtained as solutions of \eqref{eq_modified_invcone_autonomous} but can be easily computed using a linear analysis. At specific energy levels, solutions start crossing the switching hyperplane, thereby introducing the effect of the nonlinearity into the response. Hence, the response frequencies increase rapidly with the energy. For increasing energy, these frequencies tend to the limits obtained by the modes of the system with a vanishing gap, as was shown in many studies such as \cite{Jiang2004} and \cite{Moussi2015}. 
Regarding the stability properties of the NNMs, the same behavior reported in \cite{Jiang2004} is confirmed by our analysis using the Floquet multipliers obtained from the monodromy matrix defined in \eqref{eq_monodromy_expm}. As illustrated in \ref{fig:backbone_mode1}, the first NNM becomes unstable at around $\omega \approxeq 0.77$, before recovering stability at around $\omega \approxeq 0.783$. The second NNM is shown to lose stability at $\omega \approxeq 1.908$. To highlight the connection between the invariant cones and the NNMs, we note that every point on the nonlinear part of the backbone curves is obtained as a solution of the modified invariant cone problem \eqref{eq_modified_invcone_autonomous}. This corresponds to solving the classic invariant cone obtained from \eqref{eq_augmented_invcone_autonomous} for the augmented homogeneous system \eqref{eq_augmented_hom_sys_autonomous} and constraining the solution to satisfy \eqref{eq_constraining_gaplength}. 
\begin{figure}
	\centering	
	\includegraphics{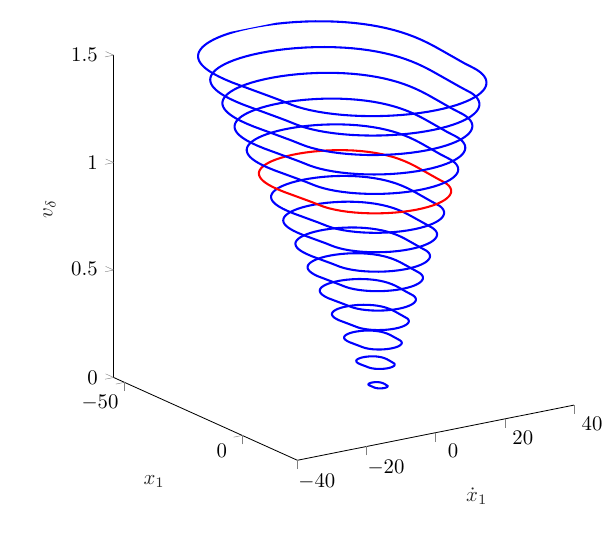}	
	\caption{Invariant cone with frequency $\omega=0.8105$. The red periodic orbit for $v_\delta=\delta=1$ corresponds to the first NNM.}
	\label{fig_3d_cones}
\end{figure}
This is visualized in Figure \ref{fig_3d_cones} for the first NNM at the frequency $\omega = 0.8105$. A projection of trajectories generated from the invariant half-lines in the augmented state space $\mR^{5}$ is shown in the three-dimensional space $(x_1,\dot{x}_1,v_\delta)$. Therein, one can see that the periodic orbits take place on planes at constant values of $v_\delta$. Scaling the initial conditions along with the gap length yields a family of mechanical systems with varying clearance gaps sharing the same response period. Obviously, increasing the gap length $v_\delta$ leads to increasing the initial values for the state vector $\vx$ to achieve the same period of the motion. In order to obtain the NNM of the actual system \eqref{eq_Jiang2004_eom}, it suffices to fix the additional dimension $v_\delta$ to its real value $\delta=1$ and recover the corresponding periodic orbit shown in red in Figure~\ref{fig_3d_cones}. This is repeated along the backbone curve for varying energy levels, which is equivalent to solving the modified invariant cone problem \eqref{eq_modified_invcone_autonomous}. This, in turn, implies solving the augmented invariant cone problem \eqref{eq_augmented_invcone_autonomous} and then picking the periodic solution corresponding to the real value of the gap length. 
\subsubsection*{Invariant manifold}
\begin{figure}
	\centering	
	\includegraphics[scale=1]{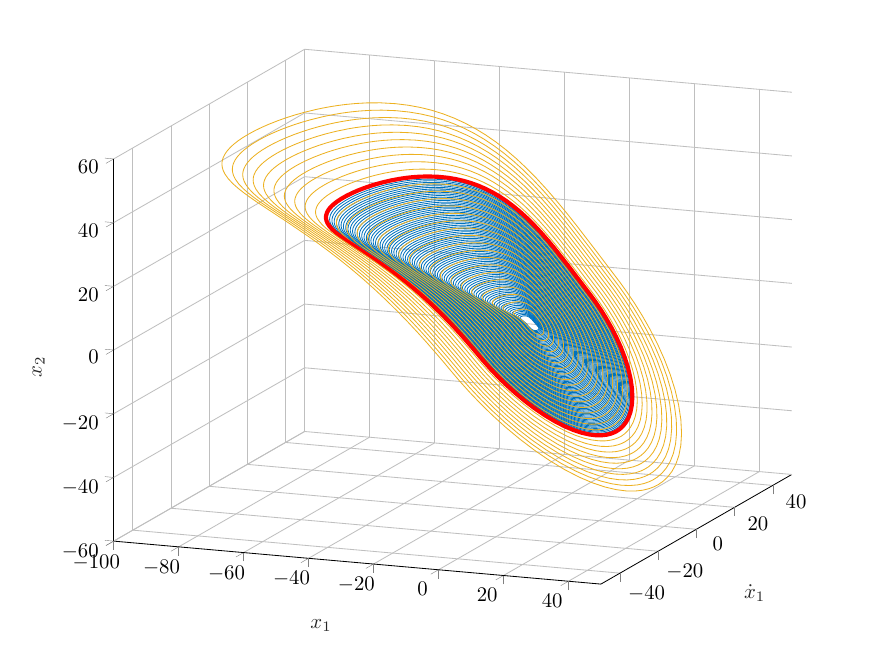}	
	\caption{Invariant manifold for the first NNM.}
	\label{fig_invmanifold_mode1}
\end{figure}

\begin{figure}
	\centering	
	\includegraphics[scale=1]{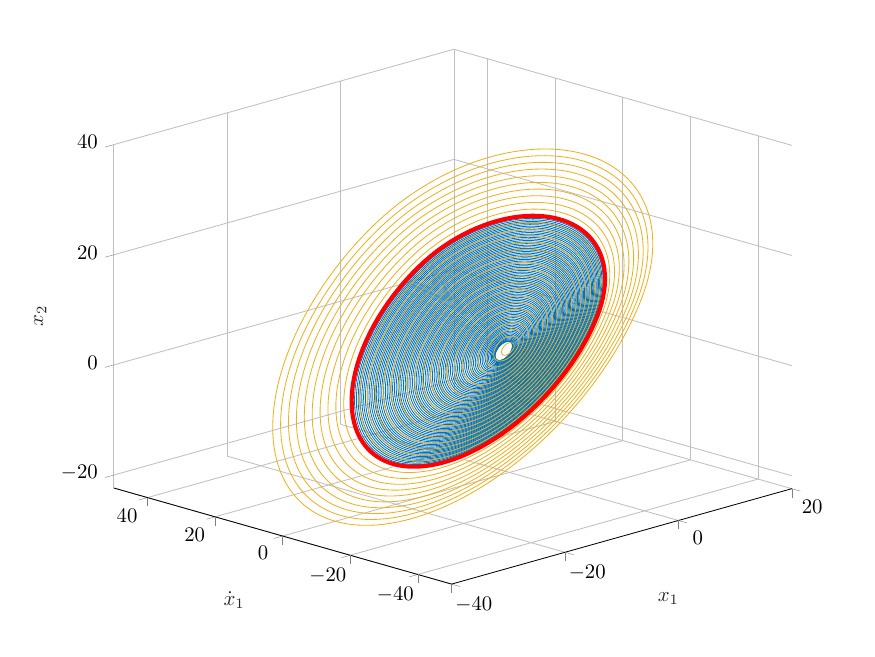}	
	\caption{Invariant manifold for the second NNM.}
	\label{fig_invmanifold_mode2}
\end{figure}

\begin{figure}
	\begin{subfigure}{0.49\linewidth}
		\centering
		\includegraphics{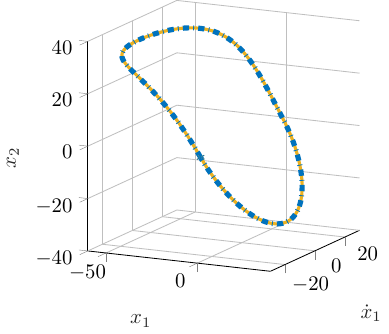}%
		\caption{Trajectories in three dimensions. }
		\label{fig:statespace_mode1_po40_cone5}
	\end{subfigure}
	\begin{subfigure}{0.49\linewidth}
		\centering
		\includegraphics{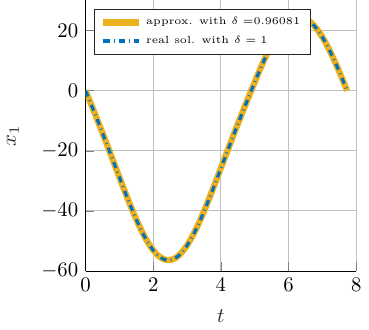}
		\caption{Time histories.}
		\label{fig:timehistory_mode1_po40_cone5}
	\end{subfigure}
	\caption{Comparing a real solution from the NNM (blue) with frequency $\omega_\textrm{real}=0.81263$ and its approximation from an invariant cone (yellow) with frequency $\omega_\textrm{cone}=0.81279$. The approximative solution was obtained by a scaling factor corresponding to $\delta=0.96081$.}
	\label{fig:vergleich_real_approx_mode1_po40_cone5}
\end{figure}

\begin{figure}
	\begin{subfigure}{0.49\linewidth}
		\centering
		\includegraphics{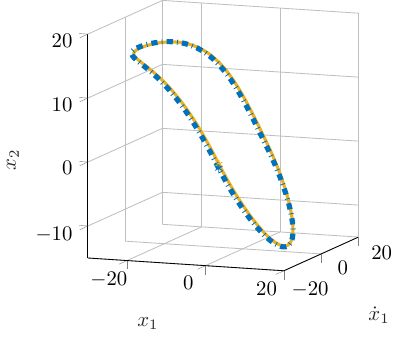}%
		\caption{Trajectories in three dimensions. }
		\label{fig:statespace_mode1_po500_cone5}
	\end{subfigure}
	\begin{subfigure}{0.49\linewidth}
		\centering
		\includegraphics{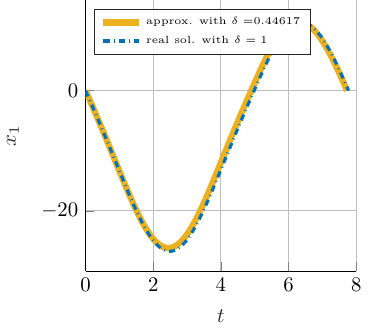}
		\caption{Time histories.}
		\label{fig:timehistory_mode1_po500_cone5}
	\end{subfigure}
	\caption{Comparing a real solution from the NNM (blue) with frequency $\omega_\textrm{real}=0.80821$ and its approximation from an invariant cone (yellow) with frequency $\omega_\textrm{cone}=0.81279$. The approximative solution was obtained by a scaling factor corresponding to $\delta=0.44617$.}
	\label{fig:vergleich_real_approx_mode1_po500_cone5}
\end{figure}

\begin{figure}
	\begin{subfigure}{0.49\linewidth}
		\centering
		\includegraphics{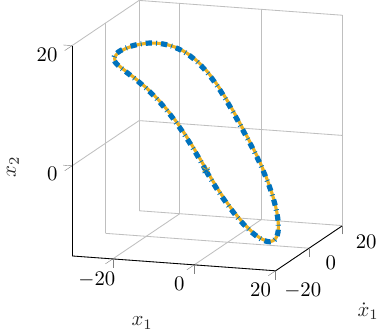}%
		\caption{Trajectories in three dimensions. }
		\label{fig:statespace_mode1_po500_cone400}
	\end{subfigure}
	\begin{subfigure}{0.49\linewidth}
		\centering
		\includegraphics{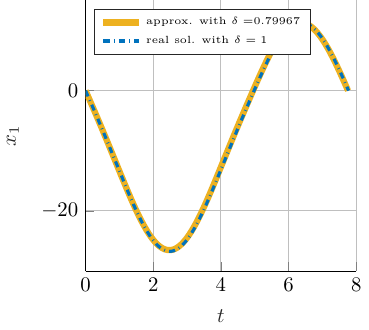}
		\caption{Time histories.}
		\label{fig:timehistory_mode1_po500_cone400}
	\end{subfigure}
	\caption{Comparing a real solution from the NNM (blue) with frequency $\omega_\textrm{real}=0.80821$ and its approximation from an invariant cone (yellow) with frequency $\omega_\textrm{cone}=0.80986$. The approximative solution was obtained by a scaling factor corresponding to $\delta=0.79967$.}
	\label{fig:vergleich_real_approx_mode1_po500_cone400}
\end{figure}

From a geometric perspective, all periodic solutions of \eqref{eq_Jiang2004_eom} live on invariant manifolds in the state space, which represent the geometries of both fundamental NNMs. Such a NNM invariant manifold is thus foliated by periodic orbits obtained as solutions of \eqref{eq_modified_invcone_autonomous}. Figures \ref{fig_invmanifold_mode1} and \ref{fig_invmanifold_mode2} show the first and second NNMs in the $(x_1,\dot{x}_1,x_2)-$space. 
The manifolds shown in blue represent the NNM invariant manifolds composed by periodic solutions of the system with $\delta=1$ having varying response periods depending on the energy level. To each of these periodic solutions belongs an invariant cone of the augmented 2CPL$_5$ system. This augmented invariant cone is foliated by periodic solutions sharing the same period and representing different 2CPL$_4$ systems with varying gap lengths $\delta$. As example, one periodic solution of the NNM is shown in red at a high energy level. Its corresponding augmented invariant cone is illustrated in yellow. The yellow orbits are periodic solutions of the augmented homogeneous system and correspond to varying gap lengths but sharing the same response period as the red curve. As explained before and as can be seen in Figure \ref{fig:backbones}, for high energy levels, the backbones tend to the bilinear limits and thus, the response frequency exhibits only a slight change with the energy. Therefore, approximating trajectories on the NNM manifolds using orbits of an augmented invariant cone describing systems with different gap lengths at a similar high energy level yields an accurate approximation, since the role of the clearance gap is insignificant for high oscillation amplitudes. Figures \ref{fig:vergleich_real_approx_mode1_po40_cone5}, \ref{fig:vergleich_real_approx_mode1_po500_cone5} and \ref{fig:vergleich_real_approx_mode1_po500_cone400} illustrate this clearly for the first NNM. In each case, a real periodic solution denoted by $\vx_\textrm{real}^\textrm{p}$ with a frequency  $\omega_\textrm{real}$ is chosen from the NNM invariant manifold shown in blue in Figure \ref{fig_invmanifold_mode1}. Its approximation is then sought using an augmented invariant cone belonging to another periodic motion $\vx_{\textrm{cone}}$ with another frequency $\omega_{\textrm{cone}}$ at a different energy level $b_\textrm{cone} = \log_{10}(\tilde{E}(\vx_{\textrm{cone}},\delta))$. Different choices of cones at different energy levels are compared to draw conclusions regarding the accuracy in approximating solutions on the NNM. In Figure \ref{fig:vergleich_real_approx_mode1_po40_cone5}, the periodic trajectory with frequency $\omega_\textrm{real} = 0.81263$ with corresponding energy level $b_\textrm{real} = \log_{10}(\tilde{E}(\vx_{\textrm{real}},\delta))=2.9639$ is chosen. An approximation $\vx_{\textrm{approx}}$ is sought from the invariant cone belonging to $\vx_{\textrm{cone}}$ with frequency $\omega_{\textrm{cone}}=0.81279$ chosen very close to that of the target solution. In order to find the specific periodic orbit $\vx_{\textrm{approx}}$ of the augmented system which lies closest to the real solution $\vx_{\textrm{real}}^{\textrm{p}}$, one has to determine the scaling factor of the augmented invariant cone such that $\vx_{\textrm{approx}}$ and $\vx_{\textrm{real}}^{\textrm{p}}$ have the same energy level. Since the energy expression \eqref{eq_energyexpression} is quadratic in the state $\vx$, this leads to the factor
\begin{equation}
	\beta_\textrm{approx} = \sqrt{\frac{10^{b_\textrm{real}}}{10^{b_\textrm{cone}}}}. 
\end{equation}  
Then, the approximative solution is found from the augmented invariant cone by scaling as follows:
\begin{equation}
	\vx_\textrm{approx} = \beta_\textrm{approx} \vx_\textrm{cone}, \quad \delta_\textrm{approx} = \beta_\textrm{approx} \, \delta. 
\end{equation}
The time histories and projections to three-dimensional representation of the state space are then compared. In Figure \ref{fig:vergleich_real_approx_mode1_po40_cone5}, one can see that the orbits are almost identical for high energy levels with an augmented cone having a similar frequency. Figure \ref{fig:vergleich_real_approx_mode1_po500_cone5} shows a comparison between a solution from the same invariant cone with $\omega_{\textrm{cone}}=0.81279$ but with a different orbit from the NNM with a frequency $\omega_\textrm{real}=0.80821$. We see that the approximation is still qualitatively good even though less accurate than in the previous case. By choosing an approximation using an invariant cone with a closer frequency $\omega_\textrm{cone}=0.80986$, the approximation accuracy is increased as shown in Figure \ref{fig:vergleich_real_approx_mode1_po500_cone400}.

\begin{figure}
	\begin{subfigure}{0.49\linewidth}
		\centering
		\includegraphics{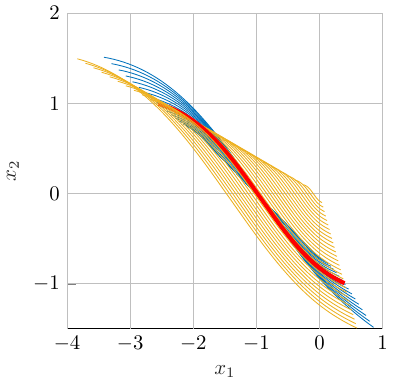}%
		\caption{First NNM and an invariant cone at $\omega= 0.6968$. }
		\label{fig:NNM_1_vs_invcone_lowenergy}
	\end{subfigure}
	\begin{subfigure}{0.49\linewidth}
		\centering
		\includegraphics{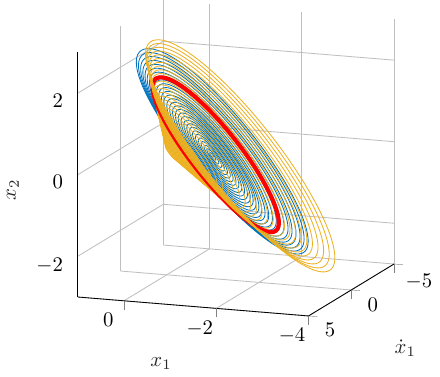}
		\caption{Second NNM and an invariant cone at $\omega = 1.8982$.}
		\label{fig:NNM_2_vs_invcone_lowenergy}
	\end{subfigure}
	\caption{Illustration of the NNM invariant manifolds (in blue) compared with an invariant cone (yellow) corresponding to a periodic solution (red) at low energy levels.}
	\label{fig:lowenergy_manifolds}
\end{figure}

\begin{figure}
	\centering
	\includegraphics{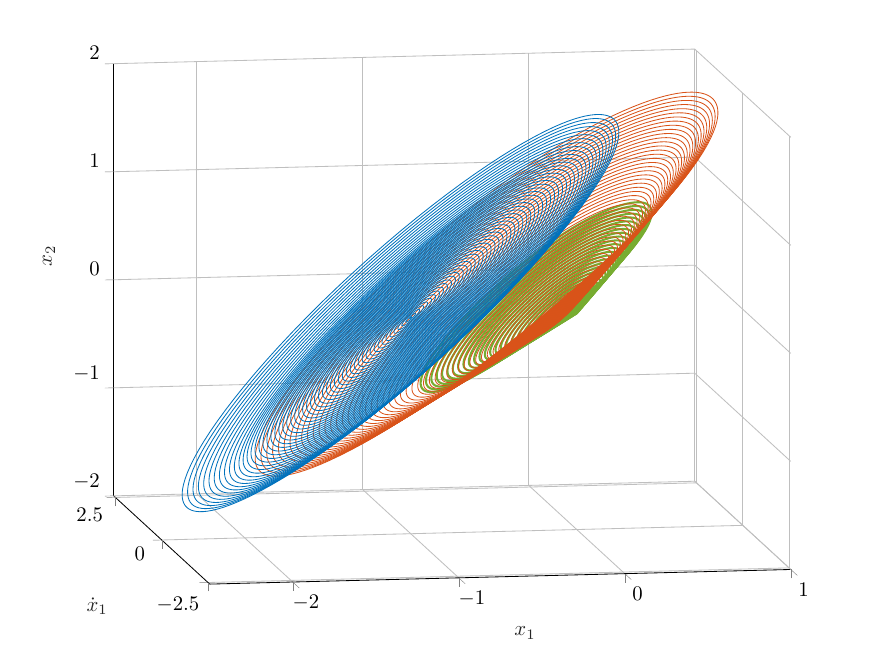}
	\caption{Zoom-in comparison of the second NNM invariant manifolds for $\delta=1$ (blue) and $\delta=0.1$ (orange) with the invariant cone of the limiting homogeneous system for $\delta=0$ (green).}
	\label{fig:NNM2_manifolds_vs_homogeneousCone}
\end{figure}

At low energy levels, the increase in the nonlinear frequency is rapid and the role played by the inhomogeneous term is more important. Therefore, only solutions with the actual gap length of the system can be picked from their corresponding augmented invariant cone. All other periodic orbits of the cone have the same period and therefore fail to accurately approximate neighboring solutions on the NNM. This is illustrated in Figure \ref{fig:lowenergy_manifolds}, where the intersection of an augmented invariant cone with the NNM manifold is pictured in red. Note that the augmented invariant cones ultimately fail to approximate the linear limit of the NNMs around the equilibrium of the underlying linear dynamics, where the NNMs are planar in contrast to the conic form of the augmented invariant cones shown in yellow. 
Decreasing the gap length leads to a smaller linear region in the NNM. The homogeneous limit is reached for a vanishing gap, which leads to the destruction of the planar region corresponding to purely linear motions on the LNM. This limiting homogeneous system, in turn, admits a classic invariant cone. The NNM of the homogeneous system in limit $\delta \downarrow 0$ degenerates to the invariant cone of the homogeneous system ($\delta=0$). Therefore, the invariant cone of the homogeneous system can be viewed as a singularity in the theory of NNMs of 2CPL$_n$ systems. 
This is illustrated in Figure \ref{fig:NNM2_manifolds_vs_homogeneousCone}, where three-dimensional zoomed views of the second NNM manifold are compared to the invariant cone of the limiting homogeneous system. Therein, the invariant manifold for $\delta=1$ is shown in blue, whereas the NNM manifold for $\delta=0.1$ is illustrated in orange. The green manifold represents the invariant cone of the limiting system for $\delta = 0$. One can clearly see that for this very low energy level, the blue manifold is mainly planar, which implies purely linear oscillations on the underlying LNM. For the smaller gap length $\delta = 0.1$, the size of the planar region of the NNM shrinks and the NNM tends toward the invariant cone of the homogeneous limiting system with a vanishing clearance gap. 
\section{Computing forced response curves using invariant cones}\label{section:frc}
The most important motivation to compute NNMs in general nonlinear systems is that forced resonances occur in their vicinity \cite{Kerschen2009}. Furthermore, the damped dynamics is influenced by the NNMs of the underlying conservative unforced system. In the following, an extension of the approach proposed in the previous section is introduced to take external harmonic forcing into account with the aim of computing period-$1$ solutions. This is accomplished through a state space augmentation, which not only considers the gap length as an additional state, but also takes the known forcing term as an additional degree-of-freedom. We hereby note that 
in the following, we consider the dynamics described by \eqref{eq_lurelike}. A new augmented state space is introduced as 
\begin{equation}\label{eq_new_ext_state_z}
	\vz\T = \begin{pmatrix}
	\vx\T & z_c & z_s & z_\delta
	\end{pmatrix} \in \mR^{n+3}, \quad \mbox{with } \quad z_c = \cos(\Omega t+\varphi_0), \quad z_s = \sin(\Omega t + \varphi_0).
\end{equation}
Here, it is important to note that the phase variable $\varphi_0$ will play an important role in the following discussion.
Using the relations $\dot{z}_c= - \Omega z_s$ and $\dot{z}_s= \Omega z_c$, the system can be put in an autonomous and homogeneous form as 
\begin{equation}\label{eq_ExtDyn_z}
	\dot{\vz} = \begin{cases}
	\tilde{\vA}_{\textrm{ext}}^+(\Omega) \vz ,\quad \mbox{for } \bar{\ve}_1\T \vz = x_1 \geq 0 \\ 
	\tilde{\vA}_{\textrm{ext}}^-(\Omega) \vz ,\quad \mbox{for } \bar{\ve}_1\T \vz = x_1 < 0 
	\end{cases},
\end{equation}
with $\bar{\ve}_j$ the $j$-th basis vector of $\mR^{n+3}$. The initial conditions have to be set to $z_c(0) = \cos(\varphi_0)$, $z_s(0)=\sin(\varphi_0)$ and $z_\delta(0)=\delta$, in order to initialize the dynamics of the system with the correct gap length and time-dependent forcing. The system matrices in \eqref{eq_ExtDyn_z} are given as 
\begin{equation}
	\tilde{\vA}_\textrm{ext}^\pm(\Omega) = \begin{pmatrix}
	\vA^\pm & \vb_t & \bm{0} & \vb \\ \bm{0}\T & 0 & -\Omega & 0 \\ \bm{0}\T & \Omega & 0 & 0 \\ \bm{0}\T & 0 & 0 & 0
	\end{pmatrix} \in \mR^{(n+3)\times(n+3)}.
\end{equation}
Note that the parameter $\alpha$ steering the damping in \eqref{eq_lurelike} is now taken as a known value, which is incorporated in the matrices $\vA^\pm$. The dependency on the excitation frequency is kept as a system parameter, in order to allow the description for varying forcing frequencies. Continuity is automatically inherited from the fact that $\vA^\pm$ satisfy the continuity condition on the switching hyperplane. In the following, we seek periodic solutions with a period time equal to the forcing period and consisting of two distinct intervals $t^-$ and $t^+$, where the system oscillates in the $\ominus$ and $\oplus$ subsystem, respectively. Therefore, it holds that $t^- + t^+ = 2\pi/\Omega$. 
System \eqref{eq_ExtDyn_z} is a 2CPL$_{n+3}$ system. Therefore, and assuming that solutions cross the switching hyperplane from $\ominus$ to $\oplus$ once per period, the Poincaré map can be expressed as
\begin{equation}
	\mathcal{P}_{\textrm{ext}}(\vz) = e^{t^+ \tilde{\vA}_\textrm{ext}^+(\Omega)} e^{t^- \tilde{\vA}_\textrm{ext}^-(\Omega)} \vz, \quad \mbox{with } \quad \vz \in \Sigma
\end{equation}
In the following, a modified invariant cone problem, similar to \eqref{eq_modified_invcone_autonomous} but taking harmonic forcing into account, is formulated in terms of the vector of unknowns $\vX\T = \begin{pmatrix}
\vz_0\T & \varphi_0 & t^- & t^+
\end{pmatrix}\in \mR^{n+6}$ as follows:
\begin{equation}\label{eq_Ext_invcone_problem}
\vF_\textrm{ext}(\vX, \Omega) = \begin{bmatrix}
\begin{bmatrix}
\vI_{n\times n} & 0 
\end{bmatrix} \left(e^{t^+ \tilde{\vA}_\textrm{ext}^+(\Omega)} e^{t^- \tilde{\vA}_\textrm{ext}^-(\Omega)} \vz_0 - \vz_0 \right) \\ 
z_{c,0} - \cos(\varphi_0) \\ 
z_{s,0} - \sin(\varphi_0) \\
z_{\delta,0} - \delta  \\ 
\bar{\ve}_1\T \vz_0 \\ 
\bar{\ve}_1\T e^{t^- \tilde{\vA}_{\textrm{ext}}(\Omega)} \vz_0 \\ 
t^- + t^+ - \frac{2\pi}{\Omega} 
\end{bmatrix} \stackrel{!}{=} \bm{0}
\end{equation}
The first equation represents the projection of the invariant cone problem $\mathcal{P}_{\textrm{ext}}(\vz_0) = \vz_0$ on the space of the original transformed coordinates $\vx$. The second and third equation represent setting the initial conditions of the additional states $z_c$ and $z_s$ to their correct values in order to recover the forcing term. This is the only condition that has to be fulfilled since periodicity of $z_{c,s}$ is automatically fulfilled in view of the last equation which fixes the period of the response to be equal to the forcing period. This shows the aforementioned important role played by the phase variable. Indeed, we seek a set of initial conditions belonging to the switching hyperplane, such that trajectories starting from the set are periodic. However, due to the explicitly time-dependent forcing term, the periodicity of trajectories starting from the invariant half-line does not hold for an arbitrary value of the phase $\varphi_0$ of the harmonic forcing at the switching hyperplane $\Sigma$. Instead, periodicity is guaranteed only in accordance with the harmonic forcing, which has to be adjusted using the phase variable $\varphi_0$, which motivates the definitions of the additional states in \eqref{eq_new_ext_state_z}. 
Regarding the fourth equation, the periodicity condition for the state $z_\delta$ is always fulfilled due to $\dot{z}_\delta=0$. Hence, only its initial value has to be fixed to the correct clearance gap length $\delta$ in order to recover the original dynamics. As discussed in Section 3, solutions on the invariant cone leave $\Sigma$ at $t_0=0$, oscillate in the $\ominus$ subsystem before returning to $\Sigma$ and crossing it at the time $t^-$. Then, solutions are governed by the $\oplus$ subsystem during a time interval $t^+$ before crossing $\Sigma$ again and closing the periodic cycle. This yields the fifth and sixth equation in \eqref{eq_Ext_invcone_problem}. The nonlinear problem \eqref{eq_Ext_invcone_problem} contains $n+6$ equations for $n+7$ unknowns including $\vX$ and the parameter $\Omega$. Therefore, a continuation algorithm can be used to follow a solution branch by regarding the forcing frequency $\Omega$ as continuation parameter. This allows to compute the period-$1$ solutions of the forced system, i.e. having the same period time as the harmonic forcing, which cross the switching hyperplane once per cycle. Motions in the purely linear regime cannot be captured due to the structure of the Poincaré map, which is not problematic since the linear response can be easily obtained by means of a transfer function if the system behaves linearly. Moreover, subharmonic responses, i.e. having periods different from the forcing period, cannot be captured by this approach.  
\subsection{Numerical example - Forced system with gap $\delta =1$}
\begin{figure}[h]
	\centering
	\includegraphics{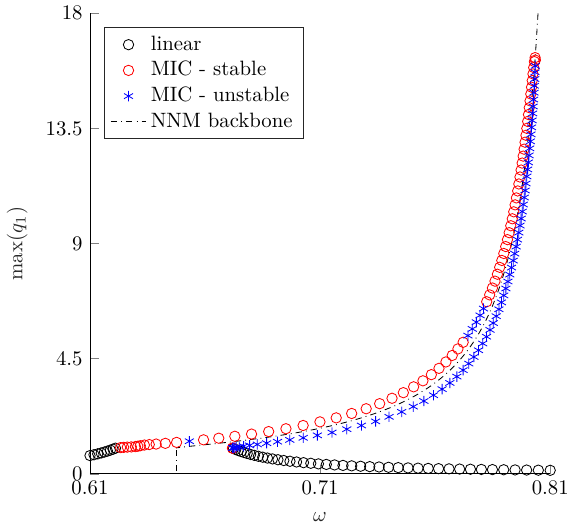}%
	\caption{Nonlinear FRC close to the first resonant frequency with forcing amplitude $0.05$.}
	\label{fig:FRC_withGap}
\end{figure}

To illustrate our approach, we revisit the example from the previous section shown in Figure \ref{fig_mech_sys_jiang2004} and add a damping term as well as external harmonic forcing applied to the first mass $m_1$. The resulting equation of motion reads as
\begin{equation}\label{eq_Jiang2004_Damped_FORCED_eom}
\vM \ddot{\vq} + \vC \dot{\vq} + \vK \vq = \vw \lambda + \vf \cos(\Omega t), \quad \mbox{with } -\lambda = \max(k_n g , 0), \quad \mbox{and } g = \vw\T \vq - \delta.
\end{equation}
The damping matrix is chosen as $\vC = 0.005 \vK$. The forcing vector is given by $\vf\T = f_\textrm{amp} \begin{pmatrix}
1 & 0
\end{pmatrix}$ with the forcing amplitude $f_\textrm{amp}=0.05$. The mass and stiffness matrix as well as the generalized force direction are kept as defined previously in \eqref{eq_Jiang_paramvalues}. The clearance gap length is also kept as $\delta = 1$. 
The same steps discussed in Section 2 can be applied to system \eqref{eq_Jiang2004_Damped_FORCED_eom} with the aim of obtaining a modified invariant cone problem of the form \eqref{eq_Ext_invcone_problem}. The forced response of the system can be determined using a nonautonomous shooting method or by means of a transfer function as long as the behavior is purely linear. The point above which solutions start crossing the switching hyperplane is then detected and used as a starting point for the continuation problem \eqref{eq_Ext_invcone_problem}. We investigate the nonlinear FRCs close to the resonant frequency of the first NNM. The nonlinear region of the FRC is obtained using the modified invariant cone (MIC) approach as solution of \eqref{eq_Ext_invcone_problem}, which yields identical results compared with the classical shooting method for nonautonomous systems. Figure \ref{fig:FRC_withGap} shows the results, where the backbone curve of the first NNM (dashed line) was superposed to the nonlinear FRC. The nonlinear FRC is generated by computing solutions of \eqref{eq_Ext_invcone_problem}, which are described in terms of transformed coordinates $\vx$ such that the switching condition is obtained by $x_1=0$. The results are then transformed back to the physical coordinates and a time integration is performed over one period time, which is also obtained from the solution as $t^- + t^+$. The maximal absolute displacement of the first mass $q_1$ is shown on the $y$-axis of the FRCs. 
One can clearly observe that the NNM backbone curve passes through the resonance peaks, which highlights once more the important role played by the NNMs in shaping the forced response of the damped system. The shape of the FRC is typical for systems with unilateral contacts, as one can see that the forced response behaves linearly up to the point where the additional stiffness is introduced in the dynamics. As a result, the resonance is shifted to higher frequencies due to the stiffening effect. Following similar arguments to the discussion in Subsection \ref{Subsection_Stability}, the stability of the obtained solutions can be directly evaluated using the monodromy matrix of the augmented system.
\subsection{Numerical example - Forced system without gap $\delta =0$}
\begin{figure}
	\begin{subfigure}{0.49\linewidth}
		\centering
		\includegraphics{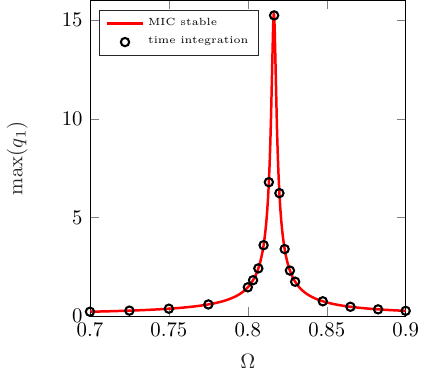}
		\caption{First resonance peak.}
		\label{fig:HomogFRC_Jiang_FirstPeak}
	\end{subfigure}
	\begin{subfigure}{0.49\linewidth}
		\centering
		\includegraphics{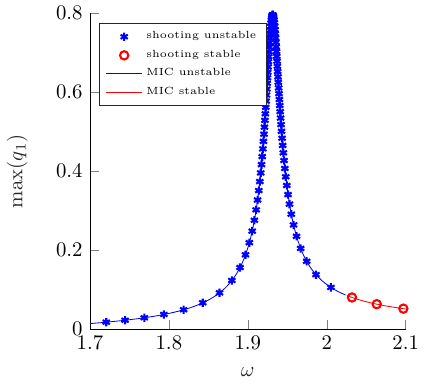}
		\caption{Second resonance peak.}
		\label{fig:HomogFRC_Jiang_SecondPeak}
	\end{subfigure}
	\caption{Nonlinear FRCs of the homogeneous system with vanishing gap $\delta=0$ with forcing amplitude $0.05$.}
	\label{fig:HomogFRC_Jiang}
\end{figure}

Here, we investigate the same 2CPL$_4$ system studied previously but set the clearance gap length to $\delta = 0$. If the state space is augmented by the additional degree-of-freedom of the harmonic forcing, the resulting augmented system becomes homogeneous. Trajectories cannot oscillate in a purely linear fashion and must therefore cross the switching hyperplane. As discussed previously, many studies have shown that the NNMs of such homogeneous systems are energy-independent, i.e. their backbone curves are straight lines with constant resonant frequencies, which is similar to linear systems. 
Since oscillations confined to only one linear subsystem are impossible when $\delta = 0$, the modified invariant cone (MIC) can capture all solutions which oscillate with the same period as the forcing, thereby covering a wider range of the nonlinear FRC. This is shown in Figure \ref{fig:HomogFRC_Jiang_FirstPeak}, where an excellent agreement is obtained between the branch of period-$1$ solutions computed using the invariant cone formulation \eqref{eq_Ext_invcone_problem} and steady-state solutions obtained using direct time integration around the first harmonic resonance. The branch of periodic solutions around the second harmonic resonance peak is shown in Figure \ref{fig:HomogFRC_Jiang_SecondPeak}, and the stability results obtained using (MIC) agree well with the shooting method. It is important to note that determining other types of periodic forced motions is beyond the scope of this work. For instance, Casini et al. \cite{Casini2012} have studied a $2$-DOF oscillator with two bilinear springs without clearance gap, i.e. a homogeneous 4CPL$_4$ system. The authors have shown that the energy-independent NNMs of the underlying unforced and undamped dynamics may bifurcate and give birth to superabundant NNMs at internal resonance conditions, which in turn can be classified into persistent or ghost modes. Persistent NNMs are those which can be retrieved in the forced response causing extra peaks in the FRCs, whereas ghost NNMs cannot be seen from the forced response even though they are stable. Other than internal resonances, forced 2CPL$_n$ systems with bilinear stiffness can exhibit sub- and superharmonic resonances, appearing as additional peaks besides the main resonances, as shown, e.g., in \cite{Karoui2023} and \cite{Kim2005}. These types of periodic motions were not considered in this work, where we mainly focused on tracing branches of fundamental NNMs and period-$1$ forced responses around the main resonances, which are characterized by a crossing behavior assumed in the structure of the Poincaré maps. 

\section{Conclusions}
In this work, we extended the concept of invariant cones, originally established for homogeneous piecewise linear systems, to take account of inhomogeneous terms in continuous PWL systems. The inhomogeneous terms can be constant or time-dependent, modeling mechanical structures with unilateral elastic contacts and under external harmonic forcing. Using an augmented state vector, the inhomogeneous dynamics are recast into an autonomous homogeneous form, thereby allowing the formulation of an augmented invariant cone problem. Special constraints are then enforced on this augmented invariant cone such that the original dynamics are recovered. The resulting modified invariant cone problem has a structure reminiscent of the classic shooting method and is suitable for tracing branches of NNMs and FRCs. One key contribution is that we established a direct link between the concept of nonlinear normal modes and the invariant cones of autonomous piecewise linear systems. This connection is threefold. First, the fundamental NNMs with one crossing per period are shown to be special solutions belonging to augmented invariant cones. These augmented invariant cones are foliated by periodic solutions of a family of mechanical systems with varying clearance gaps. The presence of a constant inhomogeneous term due to a nonvanishing gap in the original dynamics leads to NNMs characterized by a frequency-energy dependence. By taking the constant term as an additional state in the augmented dynamics, the homogeneity property governing the augmented invariant cones can be overcome to capture the interplay between the energy and the NNM frequencies of the original system. A second aspect is that the invariant cone of the limiting homogeneous system can be regarded as a singularity of the NNM. In fact, the NNMs of homogeneous PWL system have constant frequencies, as is the case for periodic solutions belonging to the corresponding invariant cone. The third facet of the link is that the invariant cones of the augmented dynamics as well as those of the limiting homogeneous case offer good approximations of the NNM invariant manifolds. This holds especially for high energy levels where the role of the inhomogeneous term becomes negligible. For low energy levels such that no crossing occurs, the NNMs are identical to the underlying planar LNMs, which cannot be captured by the invariant cones. For a decreasing inhomogeneous term, the NNM tends towards the invariant cone of the underlying homogeneous system.
Furthermore, the same methodology of augmenting the state space is also capable to take harmonic forcing into account. The resulting augmented invariant cone problem can be solved using a continuation algorithm to trace FRCs around the main resonance peaks. 
This study provides insights for the use of invariant cones to study continuous piecewise linear mechanical systems. Further work may focus on establishing a similar link between the NNMs in mechanical systems with friction and the invariant cones of discontinuous piecewise linear systems. 

\section{Conflict of interest}
The authors declare that they have no conflicts of interest.

\begin{appendices}
	\section{Shooting method and pseudo-arclength continuation}\label{Appendix_shooting}
	Based on the extended definition of Rosenberg, computing NNMs is equivalent to finding periodic orbits of the underlying undamped and unforced nonlinear system \cite{Peeters2009}. Canceling the damping and forcing terms in the equations of motion \eqref{eq_eom_2ndorder} leaves the following dynamics: 
	\begin{equation}\label{eq_cons_auto_shooting_appendix}
		\vM \ddot{\vq} + \vK \vq + \vf_{nl}(\vq) = \bm{0},  \quad \mbox{with } \vf_{nl}(\vq) = \vw \max(k_n g, 0), \quad \mbox{and } g = \vw\T \vq - \delta \in \mR, \quad \vq \in \mR^N. \nonumber
	\end{equation}
	The term $\vf_{nl}(\vq)$ contains the static nonlinear restoring force resulting from the unilateral contact. The conservative autonomous dynamics can be recast in first order form 
	\begin{equation}
		\dot{\vz} = \vg(\vz) = \begin{cases}
			\vA^+ \vz , \quad \mbox{for } g > 0 \\ 
			\vA^- \vz , \quad \mbox{for } g \leq 0,
		\end{cases} , \quad \mbox{with } \vz = \begin{bmatrix}
			\vq\T & \dot{\vq}\T 
		\end{bmatrix}\T\in \mR^{2N},
	\end{equation}
	with $\vg(\vz)$ denoting the continuous nonsmooth vector field and the matrices 
	\begin{equation}
		\vA^\pm = \begin{pNiceArray}{cc}
			\bm{0}_{N\times N} & \vI_{N\times N} \\ -\vM^{-1} \vK^\pm &  \bm{0}_{N\times N}		\end{pNiceArray}, \quad \mbox{with } \quad \vK^- = \vK ,\quad \mbox{and } \quad \vK^+ = \vK + k_n \vw \vw\T .
	\end{equation} 
	With the aim of finding periodic solutions of the governing CPL equation of motion \eqref{eq_cons_auto_shooting_appendix}, several computation methods have been established. Obviously, one may find attractive periodic solutions by simply running a numerical time integration long enough such that transient effects vanish completely and the solution converges towards a periodic motion called "steady-state". This, however, leads to a significant computational effort due to the long simulation run times and one can only find attractive periodic orbits, as unstable periodic orbits are never reached in steady-state. For this reason, periodic solution solvers have emerged as a more efficient alternative for finding periodic orbits. In this context, the shooting method is the most popular numerical technique, given its simple implementation and the ability to take nonsmooth nonlinearities into account without regularization. The underlying idea is to solve numerically the two-point boundary-value problem given by the periodicity condition 
	\begin{equation}\label{eq_shooting_equation}
		\vH(\vz_{p0},T) := \vz_p(T,\vz_{p0}) - \vz_{p0} \stackrel{!}{=} \bm{0},
	\end{equation} 
	where the function $\vH(\vz_{p0},T)$ is called the shooting function and $\vz_p(T,\vz_{p0})$ denotes the solution at time $T$ with initial condition $\vz_{p0}$.	Shooting consists in finding the initial condition $\vz_{p0}$ and the period $T$ in an iterative way, such that a periodic motion is realized. Unlike forced motion, the period $T$ of the free, i.e. autonomous, response is not known a priori. To solve the shooting condition, the method relies on both direct numerical integration and a Newton-Raphson technique. For some initial guesses $\vz_{p0}^{(0)}$ and $T^{(0)}$, the shooting equation \eqref{eq_shooting_equation} is not fulfilled and we need corrections $\Delta \vz_{p0}^{(0)}$ and $\Delta T^{(0)}$, which can be solved for using a Taylor expansion of $\vH(\vz_{p0},T)$ as 
	\begin{equation}\label{eq_taylor_shooting}
		\vH(\vz_{p0}^{(0)}),T^{(0)}) + \left. \frac{\partial \vH}{\partial \vz_{p0}}\right\rvert_{(\vz_{p0}^{(0)},T^{(0)})} \Delta  \vz_{p0}^{(0)} + \left. \frac{\partial \vH}{\partial T}\right\rvert_{(\vz_{p0}^{(0)},T^{(0)})} \Delta T^{(0)} + H.O.T = 0,
	\end{equation}
	where $H.O.T$ denotes higher order terms that are neglected. 
	This is done iteratively by replacing the superscript $(0)$ by $(k)$ and computing corrections until convergence is reached when $\vH(\vz_{p0}^{(k)},T^{(k)}) = \bm{0}$ to the desired accuracy. Since this Newton-Raphson method is a local algorithm, convergence is guaranteed only when the initial guess is sufficiently close to the solution. 
	The taylor expansion \eqref{eq_taylor_shooting} for an arbitrary iteration can equivalently be written in the form 
	\begin{equation}\label{eq_taylor_shooting_umgeformt}
		(\bm{\Phi}_T^{(k)} - I ) \Delta \vz_{p0}^{(k)} + \vg(\vz_{p0}^{(k)} ) \Delta T^{(k)} = - \vH(\vz_{p0}^{(k)}, T^{(k)} ) ,
	\end{equation}
	where $\bm{\Phi}_T$ denotes the monodromy matrix defined as 
	\begin{equation}
		\bm{\Phi}_T = \left. \frac{\partial \vz_p(t,\vz_{p0})}{\partial \vz_{p0}}\right\rvert_T.
	\end{equation}
	In smooth nonlinear systems, the monodromy matrix has to be computed from the corresponding linear ODE problem, which is best done simultaneously as running the time integration of the dynamics. In the context of 2CPL$_n$ systems, computing $\bm{\Phi}_T = \bm{\Phi}(t=T)$ reduces to solving the continuous piecewise linear ODE 
	\begin{equation}
		\dot{\bm{\Phi}} = \begin{cases}
			\vA^+ \bm{\Phi}, \quad \mbox{for } g > 0 , \\
			\vA^- \bm{\Phi}, \quad \mbox{for } g < 0 ,
		\end{cases}, \quad \mbox{with } \bm{\Phi}(0) = \vI \quad \mbox{and } \quad \vA^\pm = \begin{pNiceArray}{cc}
			\bm{0}_{N\times N} & \vI_{N\times N} \\ -\vM^{-1} \vK^\pm &  \bm{0}_{N\times N}
		\end{pNiceArray}.
	\end{equation}
	
	The right-hand side $- \vH(\bm{z}_{p0}^{(k)}, T^{(k)}) = - \bm{z}(\bm{z}_{p0}^{(k)},T^{(k)}) + \bm{z}_{p0}^{(k)}$ of equation \eqref{eq_taylor_shooting_umgeformt} is obtained by a time integration of the dynamics to yield the solution after a run time $T^{(k)}$ and comparing it with the initial condition. Therefore, each shooting iteration involves the time integration of the equations of motion to evaluate the current shooting residue. The system \eqref{eq_shooting_equation} is underdetermined since we have $2n$  equations for $2n+1$ unknowns. Therefore, the following anchor equation is usually added to remove the ambiguity due to the freedom of phase:
	\begin{equation}
		\vg(\vz_{p0}^{(k)})\T \Delta \vz_{p0}^{(k)} = 0 .
	\end{equation}
	To sum up, each Newton iteration consists in solving the linear system of $n+1$ equations 
	\begin{equation}\label{eq_shooting_updates}
		\begin{bNiceArray}{cc}
			\bm{\Phi}^{(k)} - \vI & \vg(\vz_{p0}^{(k)}) \\ \vg(\vz_{p0}^{(k)})\T & 0 
		\end{bNiceArray} \begin{bNiceArray}{c}
			\Delta \vz_{p0} \\ \Delta T
		\end{bNiceArray} =  \begin{bNiceArray}{c}
			- \bm{z}(\bm{z}_{p0}^{(k)},T^{(k)}) + \bm{z}_{p0}^{(k)} \\ 0
		\end{bNiceArray}
	\end{equation}
	for the unknown corrections $\Delta \vz_{p0}$ and $\Delta T$, after which we find the updates
	\begin{equation}
		\vz_{p0}^{(j+1)} = \vz_{p0}^{(j)} + \Delta \vz_{p0}, \quad T^{(j+1)} = T^{(j)} + \Delta T. 
	\end{equation}
	It has been already established that 2CPL$_n$ systems can exhibit at least $n$ different families of periodic orbits. In presence of a nonvanishing gap and thus an inhomogeneous term in the ODE, the NNM frequencies vary depending on the energy. A NNM family governed by the equation 
	\begin{equation}\label{eq_continuationproblem}
		\tilde{\vH}(\vz_{p0},T) = \begin{cases}
			\vH(\vz_{p0},T) = \bm{0}, \quad \mbox{shooting equation}\\
			\vg(\vz_{p0})\T \Delta \vz_{p0} = 0 , \quad \mbox{anchor equation},
		\end{cases}
	\end{equation}
	traces a curve in the $(2n+1)$-dimensional space of initial conditions $\vz_{p0}$ and periods $T$. A possible approach to follow a NNM branch is by computing a linear mode at low energy, i.e. such that the system oscillates without reaching the switching hyperplane. Then, a continuation technique can be applied to trace the branch and follow the energy-frequency curve. To this aim, one might be inclined to use the period $T$ to parametrize the NNM branch and take the previous solution as an initial guess for each step. This requires very small increments $\Delta T$ and the algorithm is not able to deal with turning points in the NNM branch. A pseudo-arclength continuation algorithm represents an efficient remedy to these problems by using a better prediction than the last computed solution. It also considers corrections to the period $T$ and thus, optimizes the path-following of the NNM branch. Starting from an initial point on the NNM branch $(\vz_{p0,(j)},T_(j))$, the pseudo-arclength continuation algorithm computes the following point on the branch in two separate steps. First, a predictor step is performed to compute a prediction $(\tilde{\vz}_{p0,(j+1)},\tilde{T}_(j+1))$ along the tangent vector to the branch as 
	\begin{equation}
	\begin{bNiceArray}{c}
	(\tilde{\vz}_{p0,(j+1)} \\ \tilde{T}_{(j+1)}
	\end{bNiceArray} = \begin{bNiceArray}{c}
	\vz_{p0,(j)} \\ T_{(j+1)}
	\end{bNiceArray} + s_{(j)} \begin{bNiceArray}{c}
	\vp_{z,(j)} \\ p_{T,(j)}
	\end{bNiceArray}
	\end{equation}
	where $s_{(j)}$ is the predictor stepsize. The tangent vector $\vp_{(j)}\T = \begin{bNiceArray}{cc}
		\vp_{z,(j)} & p_{T,(j)}
	\end{bNiceArray}$ to the branch defined by \eqref{eq_continuationproblem} is the solution of the system
	\begin{equation}
		\begin{bNiceArray}{cc}
			\bm{\Phi}^{(k)} - \vI & \vg(\vz_{p0}^{(k)}) \\ \vg(\vz_{p0}^{(k)})\T & 0 
		\end{bNiceArray} \begin{bNiceArray}{c}
			\vp_{z,(j)} \\ p_{T,(j)}
		\end{bNiceArray} = \begin{bNiceArray}{c}
		\bm{0} \\ 0
	\end{bNiceArray},
	\end{equation}
	with an additional normalization given by $||\vp_{(j)}|| = 1 $. This condition can be enforced by fixing one component of the vector $\vp_{(j)}$ and using the Moore-Penrose-Pseudoinverse to solve for the rest of its components and then divinding the resulting vector by its norm. 
	This step is global in the sense that there is a transition from one energy level to a new one, but is not iterative in contrast to the shooting iterations \eqref{eq_shooting_updates}. Following this, the second step is local and consists in performing corrections by a shooting iteration, in which the initial condition and the period are forced to be orthogonal to the predictor step. The orthogonality condition is added to the shooting procedure, which yields
	\begin{equation}
		\begin{bNiceArray}{cc}
			\bm{\Phi}^{(k)} - \vI & \vg(\vz_{p0}^{(k)}) \\ \vg(\vz_{p0}^{(k)})\T & 0 \\ 
			\vp_{z,(j)}^* & p_{T,(j)}
		\end{bNiceArray} \begin{bNiceArray}{c}
			\Delta \vz_{p0,(j+1)}^{(k)} \\ \Delta T_{(j+1)}^{(k)}
		\end{bNiceArray} =  \begin{bNiceArray}{c}
			- \vH(\vz_{p0,(j+1)}^{(k)},T_{(j+1)}^{(k)}) \\ 0 \\ 0
		\end{bNiceArray},
	\end{equation}
	where the prediction is used as initial guess, i.e. $\vz_{p0,(j+1)}^{(0)} = \tilde{\vz}_{p0,(j+1)}$ and $T^{(0)}_{(j+1)} = \tilde{T}_{(j+1)}$.
	This system is overdetermined and can be solved using a Moore-Penrose-Pseudoinverse. This corrector step is repeated until convergence is achieved, in which case the whole procedure is repeated starting from the prediction step to compute the following branch point. It should be noted, that a small prediction step size leads to too many computations, but a very large one leads to a high number of corrections or the breakdown of the Newton-Raphson iterations. Therefore, the prediction step size has to be adjusted automatically. For instance, the ratio between the number of iterations in the corrector step and the desired number of iterations can be evaluated. If it is too small, the prediction step can be increased. In case it is too large, then the prediction should be refined. The sign of the step size is chosen in order to follow the branch in the same direction. 
	
\end{appendices}

\bibliographystyle{abbrv}
\bibliography{bibliography_invco}

\begin{thebibliography}{10}

\bibitem{Attar2017}
M.~Attar, A.~Karrech, and K.~Regenauer-Lieb.
\newblock Non-linear modal analysis of structural components subjected to
  unilateral constraints.
\newblock {\em Journal of Sound and Vibration}, 389:380--410, Feb. 2017.

\bibitem{Butcher1999}
E.~A. Butcher.
\newblock Clearance effects on bilinear normal mode frequencies.
\newblock {\em Journal of Sound and Vibration}, 224(2):305--328, jul 1999.

\bibitem{Butcher2007}
E.~A. Butcher and R.~Lu.
\newblock Order reduction of structural dynamic systems with static piecewise
  linear nonlinearities.
\newblock {\em Nonlinear Dynamics}, 49(3):375--399, jan 2007.

\bibitem{carmona2005limit}
V.~Carmona, E.~Freire, E.~Ponce, J.~Ros, and F.~Torres.
\newblock Limit cycle bifurcation in {3D} continuous piecewise linear systems
  with two zones: Application to {Chua's} circuit.
\newblock {\em International Journal of Bifurcation and Chaos},
  15(10):3153--3164, 2005.

\bibitem{Carmona2002}
V.~Carmona, E.~Freire, E.~Ponce, and F.~Torres.
\newblock On simplifying and classifying piecewise-linear systems.
\newblock {\em IEEE Transactions on Circuits and Systems I: Fundamental Theory
  and Applications}, 49(5):609--620, May 2002.

\bibitem{carmona2005bifurcation}
V.~Carmona, E.~Freire, E.~Ponce, and F.~Torres.
\newblock Bifurcation of invariant cones in piecewise linear homogeneous
  systems.
\newblock {\em International Journal of Bifurcation and Chaos},
  15(08):2469--2484, 2005.

\bibitem{Carmona2006}
V.~Carmona, E.~Freire, E.~Ponce, and F.~Torres.
\newblock The continuous matching of two stable linear systems can be unstable.
\newblock {\em Discrete and Continuous Dynamical Systems}, 16(3):689--703,
  2006.

\bibitem{Casini2012}
P.~Casini, O.~Giannini, and F.~Vestroni.
\newblock Persistent and ghost nonlinear normal modes in the forced response of
  non-smooth systems.
\newblock {\em Physica D: Nonlinear Phenomena}, 241(22):2058--2067, nov 2012.

\bibitem{Casini2011}
P.~Casini and F.~Vestroni.
\newblock Characterization of bifurcating non-linear normal modes in piecewise
  linear mechanical systems.
\newblock {\em International Journal of Non-Linear Mechanics}, 46(1):142--150,
  jan 2011.

\bibitem{Caughey1990}
T.~Caughey, A.~Vakakis, and J.~Sivo.
\newblock Analytical study of similar normal modes and their bifurcations in a
  class of strongly non-linear systems.
\newblock {\em International Journal of Non-Linear Mechanics}, 25(5):521--533,
  Jan. 1990.

\bibitem{Chati1997}
M.~Chati, R.~Rand, and S.~Mukherjee.
\newblock Modal analysis of a cracked beam.
\newblock {\em Journal of Sound and Vibration}, 207(2):249--270, oct 1997.

\bibitem{Chen1996}
S.-L. Chen and S.~W. Shaw.
\newblock Normal modes for piecewise linear vibratory systems.
\newblock {\em Nonlinear Dynamics}, 10(2):135--164, jun 1996.

\bibitem{Cochelin2014}
B.~Cochelin.
\newblock Numerical computation of nonlinear normal modes using {HBM} and
  {ANM}.
\newblock In {\em Modal Analysis of Nonlinear Mechanical Systems}, pages
  251--292. Springer Vienna, 2014.

\bibitem{Hosham2016}
H.~A. Hosham.
\newblock Bifurcation of periodic orbits in discontinuous systems.
\newblock {\em Nonlinear Dynamics}, 87(1):135--148, aug 2016.

\bibitem{Jiang2004}
D.~Jiang, C.~Pierre, and S.~Shaw.
\newblock Large-amplitude non-linear normal modes of piecewise linear systems.
\newblock {\em Journal of Sound and Vibration}, 272(3-5):869--891, may 2004.

\bibitem{Jiang2005}
D.~Jiang, C.~Pierre, and S.~Shaw.
\newblock Nonlinear normal modes for vibratory systems under harmonic
  excitation.
\newblock {\em Journal of Sound and Vibration}, 288(4–5):791--812, Dec. 2005.

\bibitem{Karkar2013}
S.~Karkar, B.~Cochelin, and C.~Vergez.
\newblock A high-order, purely frequency based harmonic balance formulation for
  continuation of periodic solutions: The case of non-polynomial
  nonlinearities.
\newblock {\em Journal of Sound and Vibration}, 332(4):968--977, Feb. 2013.

\bibitem{Karoui2023}
A.~Y. Karoui and R.~I. Leine.
\newblock Model reduction of a periodically forced slow–fast continuous
  piecewise linear system.
\newblock {\em Nonlinear Dynamics}, 111(21):19703--19722, Oct. 2023.

\bibitem{Kerschen2014}
G.~Kerschen.
\newblock Computation of nonlinear normal modes through shooting and
  pseudo-arclength computation.
\newblock In {\em Modal Analysis of Nonlinear Mechanical Systems}, pages
  215--250. Springer Vienna, 2014.

\bibitem{Kerschen2009}
G.~Kerschen, M.~Peeters, J.~Golinval, and A.~Vakakis.
\newblock Nonlinear normal modes, part {I}: A useful framework for the
  structural dynamicist.
\newblock {\em Mechanical Systems and Signal Processing}, 23(1):170--194, jan
  2009.

\bibitem{Kim2005}
T.~Kim, T.~Rook, and R.~Singh.
\newblock Super- and sub-harmonic response calculations for a torsional system
  with clearance nonlinearity using the harmonic balance method.
\newblock {\em Journal of Sound and Vibration}, 281(3–5):965--993, Mar. 2005.

\bibitem{Kuepper2008}
T.~Küpper.
\newblock Invariant cones for non-smooth dynamical systems.
\newblock {\em Mathematics and Computers in Simulation}, 79(4):1396--1408,
  2008.

\bibitem{Kuepper2011}
T.~Küpper and H.~Hosham.
\newblock Reduction to invariant cones for non-smooth systems.
\newblock {\em Mathematics and Computers in Simulation}, 81(5):980--995, jan
  2011.

\bibitem{Lein2004}
R.~I. Leine and H.~Nijmeijer.
\newblock {\em Dynamics and Bifurcations of Non-Smooth Mechanical Systems},
  volume~18 of {\em Lecture Notes in Applied and Computational Mechanics}.
\newblock Springer Verlag, Berlin, 2004.

\bibitem{Manevich1972}
L.~Manevich and I.~Mikhlin.
\newblock On periodic solutions close to rectilinear normal vibration modes.
\newblock {\em Journal of Applied Mathematics and Mechanics}, 36(6):988--994,
  1972.

\bibitem{Moussi2015}
E.~Moussi, S.~Bellizzi, B.~Cochelin, and I.~Nistor.
\newblock Nonlinear normal modes of a two degrees-of-freedom piecewise linear
  system.
\newblock {\em Mechanical Systems and Signal Processing}, 64-65:266--281, dec
  2015.

\bibitem{Peeters2009}
M.~Peeters, R.~Viguié, G.~Sérandour, G.~Kerschen, and J.-C. Golinval.
\newblock Nonlinear normal modes, part {II}: Toward a practical computation
  using numerical continuation techniques.
\newblock {\em Mechanical Systems and Signal Processing}, 23(1):195--216, Jan.
  2009.

\bibitem{Rand1971a}
R.~H. Rand.
\newblock A higher order approximation for non-linear normal modes in two
  degree of freedom systems.
\newblock {\em International Journal of Non-Linear Mechanics}, 6(4):545--547,
  Aug. 1971.

\bibitem{Rand1971}
R.~H. Rand.
\newblock Nonlinear normal modes in two-degree-of-freedom systems.
\newblock {\em Journal of Applied Mechanics}, 38(2):561--561, June 1971.

\bibitem{Rand1974}
R.~H. Rand.
\newblock A direct method for non-linear normal modes.
\newblock {\em International Journal of Non-Linear Mechanics}, 9(5):363--368,
  Oct. 1974.

\bibitem{Rosenberg1960}
R.~M. Rosenberg.
\newblock Normal modes of nonlinear dual-mode systems.
\newblock {\em Journal of Applied Mechanics}, 27(2):263--268, jun 1960.

\bibitem{Rosenberg1961}
R.~M. Rosenberg.
\newblock On normal vibrations of a general class of nonlinear dual-mode
  systems.
\newblock {\em Journal of Applied Mechanics}, 28(2):275--283, June 1961.

\bibitem{Rosenberg1962}
R.~M. Rosenberg.
\newblock The normal modes of nonlinear n-degree-of-freedom systems.
\newblock {\em Journal of Applied Mechanics}, 29(1):7--14, Mar. 1962.

\bibitem{Rosenberg1964a}
R.~M. Rosenberg.
\newblock On the existence of normal mode of vibrations of nonlinear systems
  with two degrees of freedom.
\newblock {\em Quarterly of Applied Mathematics}, 22(3):217--234, 1964.

\bibitem{Rosenberg1966}
R.~M. Rosenberg.
\newblock On nonlinear vibrations of systems with many degrees of freedom.
\newblock In {\em Advances in Applied Mechanics Volume 9}, pages 155--242.
  Elsevier, 1966.

\bibitem{Rosenberg1964}
R.~M. Rosenberg and J.~K. Kuo.
\newblock Nonsimilar normal mode vibrations of nonlinear systems having two
  degrees of freedom.
\newblock {\em Journal of Applied Mechanics}, 31(2):283--290, June 1964.

\bibitem{Rossi2023}
M.~Rossi, A.~Piccolroaz, and D.~Bigoni.
\newblock Fusion of two stable elastic structures resulting in an unstable
  system.
\newblock {\em Journal of the Mechanics and Physics of Solids}, 173:105201,
  Apr. 2023.

\bibitem{Saad1992}
Y.~Saad.
\newblock Analysis of some {K}rylov subspace approximations to the matrix
  exponential operator.
\newblock {\em SIAM Journal on Numerical Analysis}, 29(1):209--228, Feb. 1992.

\bibitem{Saito2011}
A.~Saito and B.~I. Epureanu.
\newblock Bilinear modal representations for reduced-order modeling of
  localized piecewise-linear oscillators.
\newblock {\em Journal of Sound and Vibration}, 330(14):3442--3457, July 2011.

\bibitem{Shaw1991}
S.~Shaw and C.~Pierre.
\newblock Non-linear normal modes and invariant manifolds.
\newblock {\em Journal of Sound and Vibration}, 150(1):170--173, Oct. 1991.

\bibitem{Shaw1993}
S.~Shaw and C.~Pierre.
\newblock Normal modes for non-linear vibratory systems.
\newblock {\em Journal of Sound and Vibration}, 164(1):85--124, jun 1993.

\bibitem{Shaw1994}
S.~Shaw and C.~Pierre.
\newblock Normal modes of vibration for non-linear continuous systems.
\newblock {\em Journal of Sound and Vibration}, 169(3):319--347, Jan. 1994.

\bibitem{Uspensky2013}
B.~Uspensky and K.~Avramov.
\newblock Nonlinear modes of piecewise linear systems under the action of
  periodic excitation.
\newblock {\em Nonlinear Dynamics}, 76(2):1151--1156, Dec. 2013.

\bibitem{Uspensky2014}
B.~Uspensky and K.~Avramov.
\newblock On the nonlinear normal modes of free vibration of piecewise linear
  systems.
\newblock {\em Journal of Sound and Vibration}, 333(14):3252--3265, jul 2014.

\bibitem{Uspensky2019}
B.~Uspensky, K.~Avramov, and O.~Nikonov.
\newblock Nonlinear modes of piecewise linear systems forced vibrations close
  to superharmonic resonances.
\newblock {\em Proceedings of the Institution of Mechanical Engineers, Part C:
  Journal of Mechanical Engineering Science}, 233(23-24):7489--7497, aug 2019.

\bibitem{Vakakis1992}
A.~Vakakis.
\newblock Non-similar normal oscillations in a strongly non-linear discrete
  system.
\newblock {\em Journal of Sound and Vibration}, 158(2):341--361, Oct. 1992.

\bibitem{Vakakis1997}
A.~Vakakis.
\newblock Non-linear normal modes and their applications in vibration theory:
  an overview.
\newblock {\em Mechanical Systems and Signal Processing}, 11(1):3--22, Jan.
  1997.

\bibitem{Vakakis1996}
A.~F. Vakakis, L.~I. Manevitch, Y.~V. Mikhlin, V.~N. Pilipchuk, and A.~A.
  Zevin.
\newblock {\em Normal Modes and Localization in Nonlinear Systems}.
\newblock Wiley, jul 1996.

\bibitem{Vestroni2008}
F.~Vestroni, A.~Luongo, and A.~Paolone.
\newblock A perturbation method for evaluating nonlinear normal modes of a
  piecewise linear two-degrees-of-freedom system.
\newblock {\em Nonlinear Dynamics}, 54(4):379--393, feb 2008.

\bibitem{Weiss2012}
D.~Weiss, T.~Küpper, and H.~Hosham.
\newblock Invariant manifolds for nonsmooth systems.
\newblock {\em Physica D: Nonlinear Phenomena}, 241(22):1895--1902, nov 2012.

\bibitem{Zuo1994}
L.~Zuo and A.~Curnier.
\newblock Non-linear real and complex modes of conewise linear systems.
\newblock {\em Journal of Sound and Vibration}, 174(3):289--313, July 1994.

\end{thebibliography}
\end{document}